\documentclass[preprint,authoryear,10pt]{elsarticle}

\usepackage[utf8]{inputenc}
\usepackage{natbib}
\usepackage{graphicx}
\usepackage{textcomp}
\usepackage{gensymb}
\usepackage{amsmath}
\usepackage{amssymb}
\usepackage{amsfonts}
\usepackage{url}
\usepackage{fancyhdr}
\usepackage{float}
\usepackage{wrapfig}
\usepackage{enumerate}
\usepackage{siunitx}
\usepackage{booktabs}
\usepackage{yhmath}
\usepackage[margin=1 in]{geometry} %márgenes documento
\usepackage{pgfplots} %graficos calidad
\usepackage{pgfplotstable} %histogramas
\usepackage{tikz} % árbol
\pgfplotsset{compat=1.17}

\usepackage[showframe]{}%
 \usepackage[fleqn]{}%
 \usepackage{nccmath}

\usepackage[linesnumbered,ruled,vlined]{algorithm2e}
\DeclareUnicodeCharacter{2212}{-}
\usepackage{etoolbox}
\makeatletter
\patchcmd{\ps@pprintTitle}% <cmd>
  {Preprint submitted}% <search>
  {Working Paper}% <replace>
  {}{}% <succes><failure>
\makeatother
\begin{document}

\begin{frontmatter}
\title{Analyzing the impact of forecast errors in the planning of wine grape harvesting operations using a multi-stage stochastic model approach}
\author[dgr]{Alejandro Milani}
\ead{milani@uc.cl}
\author[dgr]{Alejandro Mac Cawley\corref{cor1}}
\ead{amac@uc.cl}

\cortext[cor1]{Corresponding author: amac@uc.cl}
\address[dgr]{Department of Industrial and Systems Engineering. Pontificia Universidad Cat\'olica de Chile, Av. Vicu\~na Mackenna 4860, Santiago, Chile}

\journal{arXiv}

\begin{abstract}
Forecasts and future beliefs play a critical role in the harvest labor hiring planning, especially when errors in them entails fixing previous made decisions, which can carry extra costs or losses. In this article, we study the effect that errors in the forecast/belief can have in the wine grape harvest planning process and the losses of the product. Errors are reflected in the prediction of yields and in the estimation of rain transition probabilities have on the value and losses of product. Also, using a multi-stage stochastic optimization model we can study the effect that second stage decisions have on the ability fix the planning decisions, reduce product losses and generate value. In a first step, we develop a multi-stage stochastic model which considers grape growth uncertainty given a belief in future events. The model decisions variables are: hiring, firing and maintaining harvest labor through periods, and also the harvested quantities in each period and block. Once the model defines the plan for the coming epoch, some decisions are implemented and a deviation in the forecast is revealed and the decision maker can adjust future decisions and beliefs. Results indicate that the effect of the errors in yield determination is not symmetrical; underestimations of the yields have a more significant negative effect on the objective function, while overestimation does not. Flexibility to revise hiring decisions does not make a significant difference if the yields are overestimated. The model significantly reduces losses of the better-quality grapes, since they correspond to a significant proportion of the income and account for the largest portion of income loss. Last, grapes that have an early improvement of their quality give the decision-maker an extra level of flexibility to adjust the harvesting plan.
%making a forecast error cause loss of value. This is not symmetric if the grape growth is overestimated or underestimated. Flexibility makes a significant difference reducing the amount of unharvested grape and hence, the loss of income, especially when growing rates estimations differ from reality. Grapes whose quality improves early and remains optimal for a longer period, and its maximum reachable value is higher, are recommended to treat with.

\end{abstract}

\begin{keyword}
Harvesting planning\sep Multi-stage stochastic optimization \sep Uncertainty modeling \sep Forecast Errors \sep OR in agriculture
\end{keyword}

\end{frontmatter}

\section{Introduction}

Operations planning is an important step in any activity, as it aligns resources in order to achieve the optimal economic value of production. This is particularly critical in agriculture operations, where uncertainty is always present. In fact, agricultural planners need to deal with a number of uncertain factors, such as biological and environmental, among others; which can generate significant variability and add complexity to the production planning process. 

To reduce the effect that variability has on the production planning, managers seek for information and forecasting methods which aim to reduce the uncertainty of future events. However, these models generate errors which need to be handled as the state of nature reveals. Another way to handle uncertainty is the possibility or flexibility to reassign after the uncertainty reveals itself which gives a performance advantage \citep{avanzini2021comparing}. Flexibility relates to the ability to re-allocate or re-distribute resources in the most effective way, after any uncertainty has been revealed \citep{chen2018value}. These has been modeled in stochastic programming using second-stage recourse decisions which allow the decision maker to adjust the system as new information is available. Flexibility is not only a desirable characteristic, but it is quickly becoming a requirement for survival of production-oriented companies \citep{arafa2012quantifying, barad2013flexibility, chryssolouris2013flexibility, patel2012enhancing, patel2011role, shi2003survey}

The effect of forecasts in planning is not new in the literature, \cite{ziarnetzky2018rolling} studies the performance of the planning models, using  a rolling horizon environment using a simulation model of a scaled-down wafer fab that considers forecast updates and one that does not. \cite{altendorfer2016effects} study the effect  of long-term forecast error on the optimal planned utilisation factor for a production system facing stochastic demand. However in realm of agriculture, where forecast inaccuracy is known, has not been studied. 

\cite{ferrer2008optimization} and \cite{arnaout2010optimization} presents a grape harvesting optimization model which accounts for the quality degradation as the grapes are not harvested in their optimal date. In their work, they do not account for variability in the grape growth nor stock uncertainty. In a more recent work by \cite{avanzini2021comparing}, they present a multi-farm multi-period model which consider demand, maturation, harvest, and yield risk and solves an expected value problem. They found that considering the uncertainties produces value compared to not doing it. However, in this case the authors did not account for the possibility of the planner to revise its decisions as the states of the nature reveals itself. Another way of considering the uncertainty, is the use of a multistage stochastic optimization model (MSSOM) \citep{birge2011multistage, pflug2016empirical} where the decision is made for each node of a tree of events considering its history but considering possible futures. The MSSOM approach is more complex to obtain computationally but it prescribes a tree of decisions according to the evolution of the uncertainty over time. The work by \cite{ahumada2012tactical} develops a two-stage stochastic model in which the decisions in a first stage are planting constraints and the costs associated with the planting decisions, such as labor cost and availability.

In this research we will consider the planning of harvesting operations for grapes destined to wine production, where uncertainty is given by different rate of grape growth and their transition probabilities, affecting the yields and the revenue. In this setting, the decision-maker needs to decide the amount of workers to hire, its allocation and harvesting schedule, accounting for an uncertain future grape stock. We will study the effect that forecast accuracy’s or the decision maker beliefs have on  on the economic value of production planning, depending on the grape’s quality behavior and on the level of flexibility to fix taken decisions, based on a multistage stochastic optimization model. The uncertainty will be modeled by considering two main sources: different possible grape’s yields and the allocation of transition probabilities between them over periods. The model will also consider the gapes quality over time so that if it is not harvested in the optimal ripeness period the grape losses economic value \citep{ferrer2008optimization}. We present a multistage stochastic optimization model which accounts for the variability in the grape yields over time, future stock belief and the quality behavior of the product.

The contribution of this work is threefold: first, to model the effect that errors in the planners' beliefs have on the quality of the harvesting plan; second, under which conditions these mistakes have high impact on the plan and finally, how does resource flexibility affect it. To achieve this, we will compare the results of different schedules depending on the accuracy of the future belief and the level of flexibility to fix decisions after uncertainty reveals, using our MSSOM, and determine a quantitative value for  flexibility and accuracy and their impact on the economic value of production planning. Finally, we will analyze how the quality characteristics of the grape affects the production value under these scenarios. 

%The main research questions of this research are: We know that making a mistake in harvest planning leads to increase costs and have loss income. If so, under which conditions these mistakes have high costs? Under which conditions the value of flexibility is higher? Where is more convenient to invest and have more flexibility? Which type of grape have more favorable characteristics against uncertainty?

The document is organized as follows. In Section 2, we present a literature review on how uncertainty has been incorporated in production planning and how the quality has been explored. In Section 3, we present the original optimization model, then its modifications to add uncertainty and the two grape’s quality behaviors considered in the analysis. Section 4 shows the main results, for a later discussion and conclusion in Section 5.

\section{Literature Review}

Production planning under uncertainty in agricultural systems is starting to receive increasing attention from researchers and practitioners in recent years \citep{borodin2016handling}. It has been previously studied using different approaches such as: stochastic optimization, chance constrain, robust optimization or dynamic optimization. \cite{bohle2010robust} uses a robust optimization approach for wine grape harvesting scheduling optimization problem subject to several uncertainties. \cite{moghaddam2011farm} uses a stochastic optimization model with chance constrained optimization to determine the optimal number of acres of hay a farm should harvest for their own horses’ consumption, as well as how much hay to purchase and sell to maximize the total profit of the farm. \cite{borodin2014quality} presents a stochastic optimization model for the annual harvest scheduling problem of the farmers’ entire cereal crop production at optimum maturity, using the meteorological conditions as the deciding factor that affects the harvest scheduling and progress. \cite{kennedy1988principles} has a complete book in which he looks at the applications of dynamic and stochastic dynamic programming to agriculture and natural resources. Finally, a more recent work by \cite{dowson2019multi}, presents and stochastic optimization model for a dairy farm, however there are still not many applications of multi-stage stochastic optimization model (MSSOM) in the agricultural sector. 

%Some of the applications that can be found are: \cite{kazemi2010multi} who looks at a sawmill production planning problem with uncertainty in the quality of raw materials and demand; \cite{ahumada2012tactical} develops a two-stage stochastic program to plan the production and distribution of fresh agricultural products under uncertainty; \cite{lobos2016intertemporal} determined the benefits using a stochastic modeling approach in a sawmill production environment; \cite{veliz2015stochastic} presents a harvesting and road construction decisions problem in the forestry sector in the presence of uncertainty modeled as a multi-stage problem; \cite{chen2018value} who is motivated by the problem of a seed producing company and finally, the work by \cite{varas2018assessing} who looks at the problem of a wine export focused company and its bottling planning problem facing demand uncertainty. \cite{nadal2020two} develops a two-stage stochastic model for zone delineation and crop planning under uncertainty. % 

In wine production planning, a more recent work by \cite{avanzini2021comparing} presents a MSSOM model to plan the harvest operations of wine grapes where uncertainty in weather conditions can affect the quality of grapes. They consider decisions on labor allocation and harvesting schedules, bearing in mind the uncertainty of future rain. Weather uncertainty is modeled following a Markov Chain approach, in which rain affects the quality of grapes and labor productivity. Climatic factors deteriorate grape quality over time and if they are not harvested on the optimal ripeness period. Finally, they also consider the effect in labor flexibility as the differences in ability between workers, which impact how they will cope with the effects of rain. \cite{ahumada2009application} conclude in their review, that planning models in agriculture very often fail to incorporate realistic stochastic issues in the agriculture. They go further and indicate that perhaps the reason for this lack of more realistic scenarios is the added complexity of finding solutions for the resulting models. Despite their expressive ability in modeling various real-life problems, multi-stage stochastic modeling has not been widely used in practice as they are notoriously difficult to solve. There are just few very recent examples of applications of stochastic optimization in agriculture. In the work by \cite{dowson2019multi} they formulate a stochastic optimization model of a dairy farm, \cite{flores2020stochastic} develop a framework to plan planting and harvesting schedules, and \cite{avanzini2021comparing} plan the wine grape harvesting.

%{\color{blue} Reduciria esta parte de arriba de la revision, ya que esta muy extensa en torno a planificacion en agricultura y deja poco para el foco en los errores del pronostico}

An important aspect in building a multistage stochastic model consists in developing an approximate representation of the underlying uncertainty. A common representation can be performed in the form of a scenario tree \citep{heitsch2009scenario}. The process of obtaining such reduced scenarios may vary, but generally they are based upon variance reduction techniques \citep{higle1998variance, shapiro2003monte}. \cite{lohndorf2016empirical} summarizes the state of-the art of scenarios generation of multivariate random variables for sample average approximation in: quasi-Monte Carlo methods, methods based on probability metrics, and moment matching. Still these techniques are based upon the premise that we know the underlying distribution of the events, so their objective is to define a reduced representation to make the problem tractable. 

In many of the cases we do not have complete information about how the future events distribute and we can only rely on historical data to infer the distribution. In these cases, is relevant to determine how we use and process the historical information to generate the future scenarios. We can find different approaches to this: forecasting methods, clustering methods and heuristics. In the work by \cite{dowson2020partially} they present a framework where a policy graph provides a natural means for decomposing the multistage stochastic program into a collection of subproblems, with arcs linking these subproblems representing the flow of information through time. This allows to naturally account for the information update as the states of nature reveals themselves. It also allows to solve a partially observable problem with continuous state and control variables using a stochastic dual dynamic programming (SDDP) approach. 

The use historical data for the generation of future scenarios is not new in the agricultural sector. \cite{chen2004yield}, study how the climate change influence on the distribution of future crop yields, specifically they analyze the effect that the variance have on production. \cite{murynin2013analysis}, uses image sequences of over 10 years to build and compare four yield prediction models. The models are developed through gradual addition of complexity. The initial model is based on linear regression using vegetation indices while the final model is non-linear.

%{\color{blue} Ver el tema de errores en prediccion}

The impact that forecast errors have on the overall quality of the plan and value of the objective function has been studied in the literature by analyzing the optimal learning levels. \cite{ziarnetzky2018rolling} studies the performance of the planning models, using  a rolling horizon environment using a simulation model of a scaled-down wafer fab that considers forecast updates and one that does not. \cite{he2018optimal} analyzes the value of information by maximizing an objective function represented by a nonlinear parametric belief model, while simultaneously learning the unknown parameters, by guiding a sequential experimentation process which is expensive. \cite{huang2019optimal} determines that an accurate evaluation of the expected operational cost associated with an allocation decision can be very expensive. They propose a learning policy that adaptively selects the fleet allocation to learn the underlying expected operational cost function by incorporating the value of information. Related to production planning, \cite{altendorfer2016effects} study the effect  of long-term forecast error on the optimal planned utilisation factor for a production system facing stochastic demand. Forecast errors can be transmitted from the researcher to the farmer by recommendations \citep{kolajo1988forecast}. Thus, the choice of model and the assumptions incorporated to the model may constitute a source of errors. In a warehouse environment study by \cite{sanders2009quantifying}, they found that forecast biases have a considerably greater impact on organizational cost than forecast standard deviation.

There are a number of techniques to benchmark the value generated by using a MSSOM approach. \cite{huang2009value} proposes a simple way of measuring the input of the decision process: the difference between the values of the objective functions. In their work, they present the case for capacity planning, comparing the values obtained by a multistage stochastic model with a two-stage model. \cite{escudero2007value} proposes to compare the expected result of using the solution of the deterministic mode (EEV); the wait and see solution value (WS), which corresponds to the expected value of using the optimal solution for each scenario and finally; the here-and-now solution, which corresponds to the optimal solution value to the recursion problem (RP) or MSSOM. With these results we can determine the EVPI = WS – RP, which denotes the expected value of perfect information and compares here-and-now and wait-and-see and, VSS = RP − EEV which denotes the value of the stochastic solution and compares the here-and-now and expected values approaches. 

Despite all the above, none has analyzed the relation between the quality of the scenario generation and the flexibility of the system and how does that relation affect the value of using a MSSOM approach. \cite{powell2019unified} points out in his research challenges: almost no attention has been devoted to analyzing the quality of a stochastic look-ahead model. He goes further and indicates that there is a need more research to understand the impact of the different types of errors that are introduced by the approximations used. 

%{\color{blue} Hay que habalr mas del tema de Scenario Generation, errores en los beliefs y su efecto, mecanismos para mejorar las predicciones}

\section{Problem Formulation}

In this section, we present first a stylized version of \cite{ferrer2008optimization} deterministic model. Second, we will discuss the way we add uncertainty to the model in the form of grape yield variability and its behavior over time. After that, we present the reformulated model as a multistage stochastic model which considers the mentioned uncertainty. Finally we present the analysis method to estimate the forecast accuracy's impacts on the economic value of harvest planning.

\subsection{Deterministic Grape Harvesting Problem}

We propose a wine grape harvesting model based on the previous work by \cite{ferrer2008optimization}. In this work the authors present a deterministic optimization model which minimizes the labor and machine  cost and also the quality degradation of the grapes. For this, they introduce a quality loss function which generates extra costs when harvesting deviates from the ideal date. The model determines the amount of labor and machine, its assignment to each lot and day, also the kilograms of grape to be harvested and sent to the cellar.

For our model, we assume that a farmer owns $J$ blocks that contains wine grapes, with an initial stock $s_{j,0}$, expressed in $kilograms$. The goal of the farmer or farm manager is to maximize the final profit of the wine grapes. To do that, the farmer faces a decision window of $T$ periods and on each period he/she has to make a decision regarding the the amount of work entered and dismissed, $x_t$ and $y_t$ respectively and the quantity of labor to allocate $z_{j,t}$, that implies a capacity of harvesting in that period and in that block. The net labor available in period $t$ is denoted by $m_t$ and $\beta$ is the productivity of the resource ($kilograms/period$). Grapes are harvested and sent to a winery during the same period, where $K$ is the single-period reception capacity of the winery. 

The costs are mainly about labor force. There are cost of hiring $C_e$, of termination $C_f$  and costs of keeping labor between periods $C_k$. Additionally, there is the harvesting cost $C_h$ ($\$/kilograms$) that is a productivity payment. 

Income is generated by selling the harvested grapes at a market price, which we denote by $b_j$ for the grapes of block $j$. The actual price, however, is affected by the final quality of grapes, which depends on the specific harvest time. We represent this by a quality factor, $q_{j,t}$, which is equal to 1 when $t$ is the optimal period for harvesting (which depends on grape ripeness), and decreases when $t$ differs from that optimal period. This is a similar representation to the one used by  \cite{ferrer2008optimization}, later we explain the specifics of this coefficients. In this way, the total net income at any time is $\sum_{j \in J} b_jq_{j,t}h_{j,t}$, where $h_{j,t}$ is the harvested amount from block $j$ in period $t$.

The basic deterministic model is the following:

\begin{center}
\[
\begin{array}{ll}

\min & \sum\limits_{j=1}^{J} \bigg\{ \sum\limits_{t=1}^{T}  \{ (C_{h}-b_{j} q_{j,t}) h_{j,t} + C_{e}x_{t} + C_{f}y_{t} + C_{k}m_{t} \} \bigg\} \\

s.t. & \\
& \begin{array}{lll}
m_{t} = m_{t-1} + x_{t} - y_{t} & t = 1,...,T & (d1) \\
\sum\limits_{j=1}^J z_{j,t} \leq m_{t}  & t = 1,...,T & (d2) \\
\sum\limits_{j=1}^J h_{j,t} \leq \beta m_t & t = 1,..., T & (d3) \\
h_{j,t} = \beta z_{j,t} & t = 1,..., T & (d4) \\
\sum\limits_{j=1}^J h_{j,t} \leq K &  t = 1,..., T & (d5) \\
h_{j,t} \leq S_{j,0} - \sum\limits_{l=1}^{t-1} h_{j,l} &  t = 1,..., T, j = 1,..., J & (d6) \\
h_{j,t} \geq 0 & t = 1,..., T, j = 1,..., J & (d7) \\
x_{t}, y_{t}, m_{t}, z_{t} \geq 0, \in \mathbb{Z}_{+} & t = 1,..., T & (d8) \\

\end{array}
\end{array}
\]
\end{center}

The objective function computes net income (with negative sign to be solved as a minimization problem). Expression (d1) is the manpower balance, while relation (d2) limits the number of allocated resources. Relation (d3) bounds the total harvest in terms of labor capability and relation (d4) indicates that every assigned worker produces in accordance with his/her productivity level. Relation (d5) bounds the total harvest in terms of single-period reception capacity as we assume that all grapes harvested in a given period have to be processed during the period. Constraint (d6), indicates the available stock of grapes to be harvested at a given moment, which is given by the initial stock minus what has been already harvested. Finally, constraints (d7) and (d8) establish the nature of the variables.

\subsection{Uncertainty Sources}

Uncertainty in the wine grape harvesting planning can arise from a number of factors: weather conditions, harvest yield, worker or machinery productivity or market (prices), among some of the conditions which can introduce uncertainty into the process. In this research we will concentrate on the uncertainty given by the amount of grapes available for harvest or yield, defined by the stock of grapes $s_{j,t}$ at each moment $t$, and also the weather conditions, which are represented as a transition probabilities that directly affect the grape stock or yield.

\subsubsection{Grape Yields}

The grape yield depends directly on weather, pest and land conditions during the productive period and these condition account for the main sources of uncertainty which add complexity to the harvest planning process. As grapes the harvest evolve during the harvest season, the effective available stock of grapes can also change due to climatic conditions or errors in the forecast of the effective yields. In the case of the uncertainty given by the climatic conditions, the decision maker will plan taking into account these conditions, which we will define by using scenarios $\omega$. These scenarios $\omega$ will be constructed by using a factor of the available grapes given by $\alpha_{t}^{\omega}$, if the factor is 1 the available stock will be unchanged, if the factor is above 1 the grape stock increases and a factor below 1, indicates that the stock is reduced. 

To model the uncertainty in yield due to climatic and errors in yield estimation, for each possible scenario $\omega$ we model the pre-harvest grape stock $S_{j,t}^{\omega}$ of a certain block $j$ in a given period $t$, that will be determined by a yield factor $\alpha_{t}^{\omega}$ that changes the current grape stock quantity, which is given by the previous stock quantity $S_{j,t-1}^{\omega}$ minus the amount of harvested grape $h_{j,t-1}^{\omega}$:

\begin{center}
$S_{j,t}^{\omega} = \alpha_{t}^{\omega} (S_{j,t-1}^{\omega} - h_{j,t-1}^{\omega})$
\end{center}

\subsubsection{Transition Probabilities}

The scenarios $\omega$ will be modeled using a tree structure and transition probabilities to account for the variability of weather conditions over the time horizon that the forecast is constructed. The grape stock in each scenario will be modeled by using a transition matrix of the yield factor $\alpha_{t}^{\omega}$. So the weather can reduce or increase the available grapes.

We will consider only three possible yield factors $\alpha$, hence there can be nine transition probabilities associated in order to complete the possible flow over time (Table \ref{tab:alphas}). This modeling allows us to generate a scenario tree, in which, each node corresponds to a possible future stage of a certain scenario (path) to have in consideration for the harvest decisions planning. Thanks to this, the resulting solution corresponds not only to a harvest plan but also a decision policy that allow the decision-maker to adjust the previous decision plan after uncertainty is revealed and flexibility allows. 

\begin{table}[h!]
\centering
$\begin{array}{c|SSS}
  \alpha_i    & {\alpha_1} & {\alpha_2} & {\alpha_3} \\
  \hline
  \alpha_1 & p_{11}     & p_{12}     & p_{13} \\
  \alpha_2 & p_{21}     & p_{22}     & p_{23} \\
  \alpha_3 & p_{31}     & p_{32}     & p_{33}\\
\end{array}$
\caption{Matrix of transition probabilities.}
\label{tab:alphas}
\end{table}

\begin{figure} [H]
\begin{center}
\resizebox{\linewidth}{!}{
\begin{tikzpicture}[nodes={draw, circle}, ->]
 
\node{0}
    child { node {0}
        child { node {0} 
                child { node {0} 
                    child { node {0} 
                        child { node {0} }
                        child { node {1} }
                        child { node {2} }
                    }
                    child { node {1} }
                    child { node {2} }
                    }
                child [missing]
                child [missing]
                child { node {1} 
                    child { node {3} }
                    child { node {4} }
                    child { node {5} }
                    } 
                child [missing]
                child [missing]
                child { node {2} 
                    child { node {6} }
                    child { node {7} }
                    child { node {8} }
                    }
    }
        child { node {1} }
        child { node {2} }
    }
    child [missing]
    child [missing]
    child { node {1} 
        child { node {3} }
        child { node {4} }
        child { node {5} }
    }
    child [missing]
    child [missing]
    child {node {2} 
        child { node {6} }
        child { node {7} }
        child { node {8} }
    };
\end{tikzpicture} }
\caption{Grape yield scenario tree}
\label{scenario}
\end{center}
\end{figure}
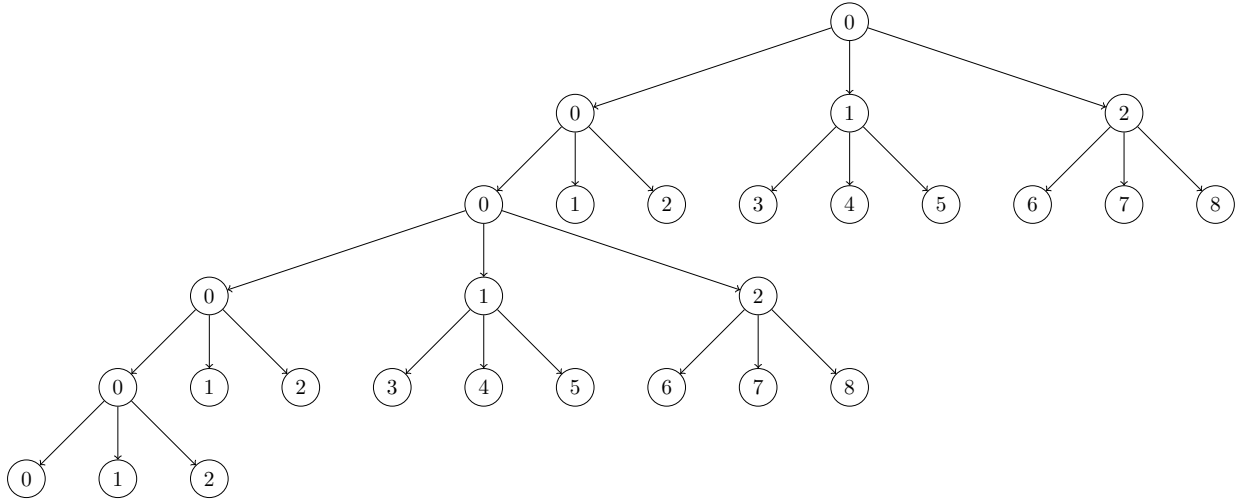

Based on the above, we estimate the probability  $p_{t}^{\omega}$ of a specific stage $\omega$ in a certain period $t$, by calculating the pitatory of the transition probabilities that shape the path from the initial node until the respective one. In Graph \ref{scenario} we can observe the representation of the aforementioned.

The occurrence probability $p_t^{\omega}$ of a certain stage $\omega$ (node) is calculated as the pitatory of the transition probabilities associated to its path from the initial stage. Defining $\mathbb{S}$ as the set of arcs $\{ij\}$ that describe the path to node n, we have:

\[p_t^{\omega} = \prod_{\{ij\} \in \mathbb{S}} p_{ij}\]

\subsection{Multi-Stage Stochastic Optimization Model}

Using the previously presented yields and transition probabilities, the multistage stochastic model that accounts for uncertainty can be modeled as follows:

\begin{itemize}
\item \textbf{Sets}
    \begin{itemize}
        \item $T$: set of periods in the time horizon.
        \item $J$: set of grape blocks of the vineyard.
        \item $\Omega$: set of future possible stages. 
    \end{itemize}
\item \textbf{Parameters}
    \begin{itemize}
        \item $b_{j}$: price of the grape in lot $j \in J$ (\$/kilograms).   
        \item $C_{e,t}$: cost of hiring in period $t \in T$ (\$/worker).
        \item $C_{f}$: cost of firing a worker (\$/worker).
        \item $C_{k}$: cost of keeping labor between periods (\$/worker per period).
        \item $C_{h}$: cost of harvesting (\$/kilograms).
        \item $S_{j,0}$: initial stock in block $j$ (kilograms). $j \in J$
        \item $q_{j,t}$: grape quality of the block $j \in J$, in the period $t \in T$.
        \item $q_{j}^{max}$: maximum reachable grape quality of the block $j \in J$.
        \item $K$: maximum reception capacity of the winery (kilograms/period).
        \item $\beta$: worker productivity (kilograms/period). 
        \item $\alpha_{t}^{\omega}$: grape yield before harvest, at time $t \in T$ in stage $\omega \in \Omega$ (positive real number).
        \item $p_{t}^{\omega}$: probability of occurrence of the stage $\omega \in \Omega$ at time $t \in T$. (positive real number $\in \{0,1\}$).
\end{itemize}
\item \textbf{Variables}
    \begin{itemize}
    \item $x_{t}^{\omega}$: workers hired at time $t \in T$ in stage $\omega \in \Omega$ (workers).
    \item $y_{t}^{\omega}$: workers laid off at time $t \in T$ in stage $\omega \in \Omega$ (workers).
    \item $m_{t}^{\omega}$: manpower or labor force at time $t \in T$ in stage $\omega \in \Omega$ (workers).
    \item $h_{j,t}^{\omega}$: harvested quantity at $j \in J$ block in period $t\in T$ in stage $\omega \in \Omega$ (kilograms/period).
    \item $S_{j,t}^{\omega}$: stock in block $j \in J$ at period $t \in T$ in stage $\omega \in \Omega$. (kilograms).
\end{itemize}

\item \textbf{Objective Function and constraints}

\begin{center}
\[
\begin{array}{ll}

\min & \sum\limits_{\omega=1}^{\Omega} \bigg\{ \sum\limits_{t=1}^{T} p_t^{\omega} \Big( C_{e,t} x_t^{\omega} + C_f y_t^{\omega}  + C_k m_t^{\omega}  + \sum\limits_{j=1}^J (C_h - b_j q_{j,t}) h_{j,t}^{\omega} \Big) +  \sum\limits_{j=1}^J p_T^{ \omega} \Big(b_j q_{j}^{max} (S_{j,T}^{\omega} - h_{j,T}^{\omega})\Big) \bigg\} \\

s.t. & \\
& \begin{array}{lll}
m_{t}^{\omega} = m_{t-1}^{\omega} + x_{t}^{\omega} - y_{t}^{\omega} & \forall t \in T, \forall \omega \in \Omega & (d1) \\
m_{t}^{\omega} = m_{t}^{\tau}, x_{t}^{\omega} = x_{t}^{\tau}, y_{t}^{\omega} = y_{t}^{\tau} & \forall t \in T, \forall \omega, \tau : \Omega^{\omega}_{[t]} = \Omega^{\tau}_{[t]} & (d2) \\
y_T^{\omega} = m_{T}^{\omega} & \forall \omega \in \Omega& (d3) \\
\sum\limits_{j=1}^J h_{j,t}^{\omega} \leq \beta m_t^{\omega} & \forall t \in T, \forall \omega \in \Omega & (d4) \\
h_{j,t}^{\omega} \leq S_{j,t}^{\omega} & \forall j \in J,\forall t \in T, \forall \omega \in \Omega & (d5) \\
\sum\limits_{j=1}^J h_{j,t}^{\omega} \leq K &  \forall t \in T, \forall \omega \in \Omega & (d6) \\
S_{j,t}^{\omega} = \alpha_{t}^{\omega} (S_{j,t-1}^{\omega} - h_{j,t-1}^{\omega}) & \forall j \in J,\forall t \in T, \forall \omega \in \Omega & (d7) \\
h_{j,t}^{\omega} \geq 0 & \forall j \in J, \forall t \in T & (d8) \\
x_{t}^{\omega}, y_{t}^{\omega}, m_{t}^{\omega} \geq 0, \in \mathbb{Z}_{+} & \forall t \in T, \forall \omega \in \Omega & (d9) \\

\end{array}
\end{array}
\]
\end{center}

\end{itemize}

The objective function computes the expected net income under uncertainty, which also considers an expected penalty for leaving unharvested grape in the last period determined by its maximum reachable quality. Expression (d1) is the manpower balance for each possible stage $\omega$, while (d2) correspond to non-anticipativity constraints. (d3) guaranties to fire the remaining harvest labor at the end of the planning. The expression (d4) bounds the total harvest in terms of labor capability and relation, while expression (d5) establishes that harvest is bounded by the volume available in the block. Relation (d6) bounds the total harvest in terms of single-period reception capacity as we assume that all grapes harvested in a given period have to be processed during the period. Expression (d7) updates the grape stock of the block from the previous period stock and harvested amount values and the yield of the respective stage. Finally, relations (d8) and (d9) establish the nature of the variables.

This allows not also make a decision plan but a harvest decision policy that indicates the best decisions subject to the actual stage conditions and possible future stages.

\section{Beliefs Error Determination}

%{\color{blue} No esta claro cómo se generan los errores en los beliefs... seria bueno hacer un subcapitulo de esto ya que es importante en el analisis}

\cite{avanzini2021comparing} studied the benefits of using a MSSOM approach when planing grape harvesting operations, however none has analyzed the effect that deviations in the grape yields forecast and transition probabilities beliefs, can have on the value of using a MSSOM optimization approach. To study the effect of such errors in the value and grape losses, in the case of grape yields, we will suppose that the decision-maker makes decisions having available just partial information, which is presented in the form of historical data, tendencies, forecasts, sensorization and so on. When the state of the nature reveals itself, the decision make will realize its error and perform actions in order to correct them, in the form of second stage decisions. The error in the grape yields will be defined as a difference in the $\alpha$ factor of grape yields values that differ from the real ones.

In the case of transition probabilities, the same phenomenon occurs when the decision-maker does not have perfect information. The believed probabilities of the scenarios are usually different from the real ones. Thus, these kind of error will be represented defining varied transition probabilities matrix that have differences between the defined as the real transition probabilities matrix.

As these two types of errors can be overestimations or subestimations of the reality, we will take into account believed values that represent that. In the case of grape yields, we will use smaller and greater values compared to the defined as real ones. In the case of transition probabilities, we will consider matrix with higher and lower probabilities values compared to the real ones.

\section{Analysis Methodology}

To determine  the effect that differences in the beliefs have on the value of using a MSSOM approach, we will analyze the effect that they have on the production performance and economic value, considering the grape quality behavior and also the level of flexibility of the decision-maker to fix the labor hiring and harvest plan.

To achieve this, we separate the analysis in four parts. The first two, studies the impacts of beliefs in grape yields, while the last two analyze the impacts of beliefs in transition probabilities. Each pair of analysis considers separately two kinds of behaviors of grape quality over time: having grapes with different quality improvement velocities, or having grapes with different maximum reachable qualities.

For the grape growth analysis, we first run the reformulated model with the real yields and save the harvest plan decisions, incomes, and costs. To this, we call it the perfect information scenario ($\mathcal{PI}$). Then, we run again the model with the believed yields and save the harvest plan decisions defined as the believed scenario ($\mathcal{BS}$). Finally, we run the model with the real yields but with the decision plan of the believed scenario, and then we save the final production values defined as the real scenario ($\mathcal{RS}$). Having these results, we compare the PI values with the RS values obtaining the economic impact of forecast accuracy.

To the previous method, we add the possibility to the decision-maker to fix the harvest plan after a certain number of periods occur and so, uncertainty reveals which allows the decision-maker to know the real yields values. The number of periods passed before being able to change the schedule has the objective to represent the level of flexibility the decision-maker has. 

For the transition probabilities analysis, we execute the same algorithm but comparing the objective function value resulting from planning based on the believed transition probabilities, with the values obtained from the real transition probabilities. To this, we also add the possibility to the decision-maker to adjust the hiring and harvest decision after a certain number of periods determined by the flexibility he/she has, knowing the real transition probabilities values.

To compare the aforementioned solutions and determine where the value is lost or gained, we will look into different components of the objective function: the objective function value itself ($\mathbb{OF}$), the objective function deficit compared to the PI scenario ($\mathbb{OFD}$), the total incomes deficit compared to the PI scenario ($\mathbb{ID}$), the percentage impact of unharvested grape ($\mathbb{UG}$) and average quality of harvested grape on total incomes ($\mathbb{AQ}$). In addition to the global components, we will analyze individually by grape block, which have different quality behaviors: the percentage variation of harvested grape compared to the PI scenario and the variation of revenue compared to the PI scenario. This will allow us to determine what type of grape should the managers seek for. 

%{\color{blue} Seria bueno el colocar pseudocodigo aca para explicitar mejor como funciona el proceso. Aca te dejo un ejemplo. O tambien un diagra que explicite el proceso iterativo.}

\begin{algorithm}[H]
 \caption{Pseudo code for optimization methodology}
 \SetAlgoLined
 \textbf{Input:} real values of yields and transition probabilities ($\alpha^{PI}_{1},\alpha^{PI}_{2},\alpha^{PI}_{3}, \mathbb{P}^{PI}$)\;
 \textbf{Run} the MSSOM\; 
    \hspace{3mm} $\mathcal{D}^{PI} \leftarrow (x_1,...,x_T,y_1,...,y_T,h_1,...,h_T)$ \;
    \hspace{3mm} $\mathbb{OF}^{PI} \leftarrow$ objective function value \;
    \hspace{3mm} $\mathbb{I}^{PI} \leftarrow$  net incomes \;
\textbf{Input:} believed values of yields or transition probabilities (according to the analysis) \;
 \textbf{Run} the MSSOM\; 
    \hspace{3mm} $\mathcal{D}^{BS} \leftarrow (x_1,...,x_T,y_1,...,y_T,h_1,...,h_T)$ \;
    \hspace{3mm} $\mathbb{OF}^{BS} \leftarrow$ objective function value \;
    \hspace{3mm} $\mathbb{I}^{BS} \leftarrow$  net incomes \; 
\textbf{Input:} ($\alpha^{PI}_{1},\alpha^{PI}_{2},\alpha^{PI}_{3}, \mathbb{P}^{PI}$), $\mathcal{D}^{BS}$, number of periods passed $f$ before re-optimize the decision plan\;
\textbf{While} $t <= f$ \textbf{do}\;
    \hspace{3mm} $x_t \leftarrow x_t^{BS} , y_t \leftarrow y_t^{BS}, h_t \leftarrow h_t^{BS}$ \;
\textbf{Run} the MSSOM for $t > f$\; 
    \hspace{3mm} $\mathcal{D}^{RS} \leftarrow (x_1^{BS},...,x_f^{BS},x_{f+1},...,x_T,y_1^{BS},...,y_f^{BS},y_{f+1},...,y_T,h_1^{BS},...,h_f^{BS},h_{f+1},...,h_T)$ \;
    \hspace{3mm} $\mathbb{OF}^{RS} \leftarrow$ objective function value \;
    \hspace{3mm} $\mathbb{I}^{RS}  \leftarrow$  net incomes \;     
\textbf{Calculate:} $\mathbb{OFD}, \mathbb{ID}, \mathbb{UG}, \mathbb{AQ}$
\end{algorithm}

\begin{fleqn}

\begin{equation}
\mathbb{OF} = 
\sum\limits_{\omega=1}^{\Omega} \bigg\{ \sum\limits_{t=1}^{T} p_t^{\omega} \Big( C_{e,t} x_t^{\omega} + C_f y_t^{\omega}  + C_k m_t^{\omega}  + \sum\limits_{j=1}^J (C_h - b_j q_{j,t}) h_{j,t}^{\omega} \Big) +  \sum\limits_{j=1}^J p_T^{ \omega} \Big(b_j q_{j}^{max} (S_{j,T}^{\omega} - h_{j,T}^{\omega})\Big) \bigg\}
\end{equation}

\begin{equation}
\mathbb{OFD} = \mathbb{OF}^{PI} - \mathbb{OF}^{RS}
\end{equation}

\begin{equation}
\mathbb{I} = 
\sum\limits_{\omega=1}^{\Omega} \sum\limits_{t=1}^{T} p_t^{\omega} \Big( \sum\limits_{j=1}^J ( - b_j q_{j,t}) h_{j,t}^{\omega} \Big) \\
\end{equation}

\begin{equation}
\mathbb{ID} = \mathbb{I}^{PI} - \mathbb{I}^{RS}
\end{equation}

\begin{equation}
\mathbb{UG} = -\mathbb{ID}/\mathbb{OFD}
\end{equation}

\begin{equation}
\mathbb{H} = \sum\limits_{\omega=1}^{\Omega} \sum\limits_{t=1}^{T} p_t^{\omega} \sum\limits_{j=1}^J h_{j,t}^{\omega}
\end{equation}

\begin{equation}
\mathbb{AQ} = \mathbb{I}^{RS} - \mathbb{I}^{PI} (\mathbb{H}^{RS}/\mathbb{H}^{PI})
\end{equation}

\end{fleqn}

\section{Model Parameters}

Below are presented the base parameters used in this work, corresponding mainly to the vineyard characteristics and planning time horizon .

\begin{table}[htbp]
  \centering
{
    \begin{tabular}{cccccc}
    \hline
    Model Parameter & Notation  &    & Value  & Units \\
    \hline
    Grape price           & $b_{j}$  &   &215     & \$/kilograms \\
    Initial harvest stock & $S_{j,0}$&   & 7,000  & kilograms \\
    Worker productivity   & $\beta$  &   & 1,600  & kilograms/period  \\
    Vineyard blocks       & $n(J)$   &   & 6      & blocks \\
    Planning time horizon & $ n(T)$  &   & 6      & periods \\
        \hline
    \end{tabular}
    }
  \caption{Model base parameters}
  \label{tab:parameters1}
\end{table}

\subsection{Costs}

We decided to determine different hiring cost in each period by using a time dependent exponential function to represent the decreasing of labor availability the later the model tries to hire and also the market price increment due to the previous argument. Therefore the mentioned function is the following:

\begin{center}
$C_{e,t} = C_e(t) = 100,000 \cdot 1.7^{t} \quad  \forall t \in T \quad  (\$/worker)$
\end{center}

Next are presented the other costs parameters considered in the model, which are based on Chilean grape market values.

\begin{table}[htbp]
  \centering
{
    \begin{tabular}{cccccc}
    \hline
    Model Parameter & Notation  &    & Value  & Units \\
    \hline
    Cost of firing & $C_{f}$ &       & 4,000   & \$/worker\\
    Cost of keeping labor &$C_{k}$&  & 3,500   & \$/worker \\
    Cost of harvesting & $C_{h}$ &   & 21      & \$/kilograms  \\
        \hline
    \end{tabular}
    }
  \caption{Model costs parameters}
  \label{tab:parameters2}
\end{table}

\subsection{Grape Yields}

The model builds the harvest schedule considering three possible grape yield scenario. As the decision-maker does not have perfect information of the future, the believed grape yields almost always do not match with the real values. To represent these value errors, we consider weights that relate the believed rates with the real ones. These is performed by changing the estimated value $\alpha_{t}^{\omega}$ in ten yield function. 

Looking at the wine grape literature, we can determine that under optimal weather conditions the availability of grapes will increase by twofold, hence the value of $\alpha_{t}^{\omega}$ is 1.5; under regular weather conditions, the availability of grapes will stay the same during the harvest season,; finally, under bad weather conditions the availability will be 0.7 times. These represent the volume increasing rate within a single period of time $t$, where each period represents a hole month. Then, the real possible yields are the same, weighted by a positive real number. When this is greater than 1, means the belief underestimated the possible yields. By the other hand, if this weight is less than 1 means the forecast overestimated the possible yields. The weights we will use in these analysis are: $\{3,2,1.5,1.1,1,0.9,0.7,0.5\}$

\subsection{Transition Probability scenarios}

Transition probabilities determine how the grape growth rate will vary over periods, since they determine if yields will be increased or reduced in the next period and by which factor ($\alpha_i$). These probabilities are forecasted by the producer, as the grape yields, and can also differ from what they were believed to be. 

For our base case, we will suppose that each growth rate has equally probability of occurrence (Scenario 4) and this are the current beliefs of the producer. From this base scenario, the real transition probabilities can differ, ranging from the pessimist scenario (Scenario 1) where initially the probabilities were thought to be equi-probable and finally the growth factor is absorbed to a value of 2. To the other extreme, the optimist scenario (scenario 8) were again it was thought to be equi-probable, but finally the growth factor is absorbed to a value of 1.1. Between these two extreme scenarios there is a continuum of scenarios of different transition probabilities which represent different levels of pessimism scenarios. 

\begin{center}
    \(
    \begin{array}{c|SSS}
                     \multicolumn{4}{c}{\text{Scenario 1}}\\
         \multicolumn{4}{c}{\text{(Optimist)}}\\
          \alpha_i & {0.7} & {1} & {1.5} \\
          \hline
          0.7 & 1   & 0     & 0 \\
          1 & 1   & 0     & 0 \\
          1.5   & 1   & 0     & 0 \\
    \end{array}
    \)
    \hspace{.2in}
    \(
    \begin{array}{c|SSS}
                     \multicolumn{4}{c}{\text{Scenario 2}}\\
         \multicolumn{4}{c}{}\\
          \alpha_i & {0.7} & {1} & {1.5}\\
          \hline
          0.7 & 0.5   & 0.5     & 0 \\
          1 & 0.5   & 0.5     & 0 \\
          1.5   & 0.5   & 0.5     & 0 \\    
    \end{array}
    \)
     \hspace{.2in}
    \(
    \begin{array}{c|SSS}
                     \multicolumn{4}{c}{\text{Scenario 3}}\\
         \multicolumn{4}{c}{}\\
          \alpha_i & {0.7} & {1} & {1.5} \\
          \hline
          0.7 & 0.6   & 0.3     & 0.1 \\
          1 & 0.4   & 0.3     & 0.3 \\
          1.5   & 0.5   & 0.3     & 0.2 \\
    \end{array}
    \)
\end{center}
\begin{center}
    \(
    \begin{array}{c|SSS}
                     \multicolumn{4}{c}{\text{Scenario 4}}\\
         \multicolumn{4}{c}{}\\
          \alpha_i & {0.7} & {1} & {1.5}\\
          \hline
          0.7 & 0.3   & 0.6     & 0.1 \\
          1 & 0.3   & 0.4     & 0.3 \\
          1.5   & 0.3   & 0.5     & 0.2 \\
    \end{array}
    \)
    \hspace{.2in}
    \(
    \begin{array}{c|SSS}
                     \multicolumn{4}{c}{\text{Scenario 5}}\\
         \multicolumn{4}{c}{\text{Base case}}\\
          \alpha_i & {0.7} & {1} & {1.5} \\
          \hline
          0.7   & {0.}\overline{3} & {0.}\overline{3} & {0.}\overline{3}  \\
          1 & {0.}\overline{3} & {0.}\overline{3} & {0.}\overline{3}  \\
          1.5 & {0.}\overline{3} & {0.}\overline{3} & {0.}\overline{3}  \\
    
    \end{array}
    \)
     \hspace{.2in}
    \(
    \begin{array}{c|SSS}
                 \multicolumn{4}{c}{\text{Scenario 6}}\\
         \multicolumn{4}{c}{}\\
          \alpha_i & {0.7} & {1} & {1.5}\\
          \hline
          0.7 & 0.1   & 0.3     & 0.6 \\
          1 & 0.3   & 0.3     & 0.4 \\
          1.5   & 0.2   & 0.3     & 0.5 \\
    \end{array}
    \)
\end{center}
\begin{center}
    \(
    \begin{array}{c|SSS}
             \multicolumn{4}{c}{\text{Scenario 7}}\\
         \multicolumn{4}{c}{}\\
          \alpha_i & {0.7} & {1} & {1.5}\\
          \hline
          0.7 & 0   & 0.5     & 0.5 \\
          1 & 0   & 0.5     & 0.5 \\
          1.5   & 0   & 0.5     & 0.5 \\
    \end{array}
    \)
    \hspace{.2in}
    \(
   \begin{array}{c|SSS}
         \multicolumn{4}{c}{\text{Scenario 8}}\\
         \multicolumn{4}{c}{\text{(Pessimist)}}\\
          \alpha_i & {0.7} & {1} & {1.5} \\
          \hline
          0.7 & 0   & 0     & 1 \\
          1 & 0   & 0     & 1 \\
          1.5   & 0   & 0     & 1 \\
    
    \end{array}    
    \)
\end{center}

\subsection{Grape Quality Characteristics}

%{\color{blue} Esta algo confuso la diferencia entre ambos graficos, ya que tienen mismos periodos pero distntas fechas en que se llega al optimo. Seria conveniente detallar ambos}

Gape quality has an effect on the harvest decisions as previously studies have pointed out. To study this factor we will divide the quality aspect into two characteristics. First, we will look at different maximum levels of quality achievable by the grapes and how this affects the harvest decisions.In a second aspect, we will keep the maximum quality level fixed, but we will vary the speed in which the grape reaches the optimal quality conditions.

As we know, quality is a factor that affects the market value of the grape, which is why it was modeled as a weight to that price whose value varies between zero and one. The closer to one this is, the better quality grapes are.

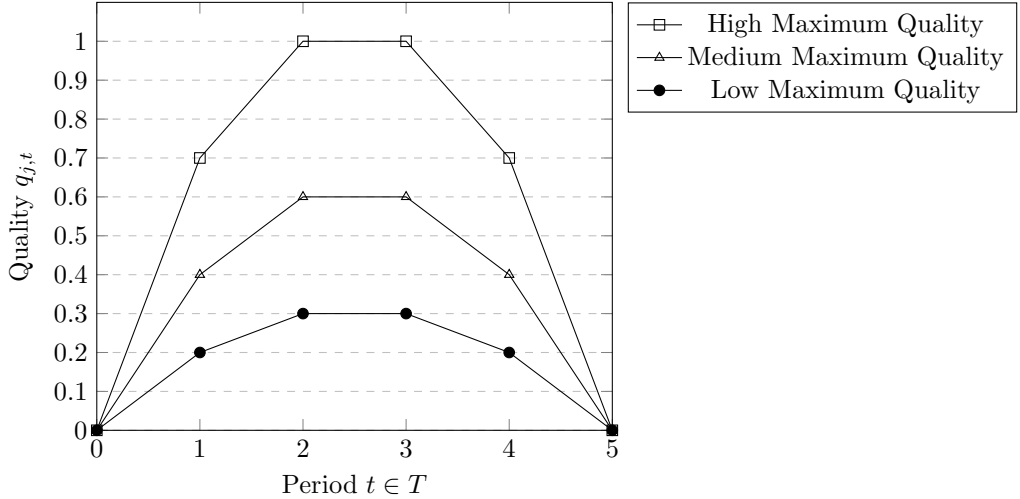
\begin{figure} [H]
\begin{center}
\begin{tikzpicture}
\begin{axis}[
    xlabel={Period $t \in T$},
    ylabel={Quality $q_{j,t}$},
    xmin=0, xmax=5,
    ymin=0, ymax=1.1,
    xtick={0,1,2,3,4,5},
    ytick={0,0.1,0.2,0.3,0.4,0.5,0.6,0.7,0.8,0.9,1},
    legend pos = outer north east,
    ymajorgrids=true,
    grid style=dashed, ]

\addplot[color=black,mark=square,    ]
    coordinates {(0,0)(1,0.7)(2,1)(3,1)(4,0.7)(5,0) };

\addplot[color=black,mark=triangle,    ]
    coordinates {(0,0)(1,0.4)(2,0.6)(3,0.6)(4,0.4)(5,0) };
    
\addplot[color=black,mark=*,    ]
    coordinates {(0,0)(1,0.2)(2,0.3)(3,0.3)(4,0.2)(5,0) };
    
\legend{High Maximum Quality, Medium Maximum Quality, Low Maximum Quality},

\end{axis}
\end{tikzpicture}
\caption{Grape quality behavior: Different maximum reachable qualities by block}
\label{q1}
\end{center}
\end{figure}

In line with the objectives mentioned above, for the first analysis to simulate having different maximum reachable qualities we considered grape blocks whose qualities all start and end in zero, but having different maximum reachable values. Graph \ref{q1} shows a representation of this case.

\begin{figure} [H]
\begin{center}
\begin{tikzpicture}
\begin{axis}[
    xlabel={Period $t \in T$},
    ylabel={Quality $q_{j,t}$},
    xmin=0, xmax=5,
    ymin=0, ymax=1.1,
    xtick={0,1,2,3,4,5},
    ytick={0,0.1,0.2,0.3,0.4,0.5,0.6,0.7,0.8,0.9,1},
    legend pos = outer north east,
    ymajorgrids=true,
    grid style=dashed, ]

\addplot[color=black,mark=square,    ]
    coordinates {(0,0)(1,0.6)(2,0.9)(3,1)(4,1)(5,1) };

\addplot[color=black,mark=triangle,    ]
    coordinates {(0,0)(1,0.2)(2,0.4)(3,0.6)(4,0.8)(5,1) };
    
\addplot[color=black,mark=*,    ]
    coordinates {(0,0)(1,0)(2,0)(3,0.1)(4,0.4)(5,1) };
    
\legend{High Ripening Rate, Medium Ripening Rate, Low Ripening Rate},

\end{axis}
\end{tikzpicture}
\caption{Grape quality behavior: Different ripening rates by block}
\label{q2}
\end{center}
\end{figure}
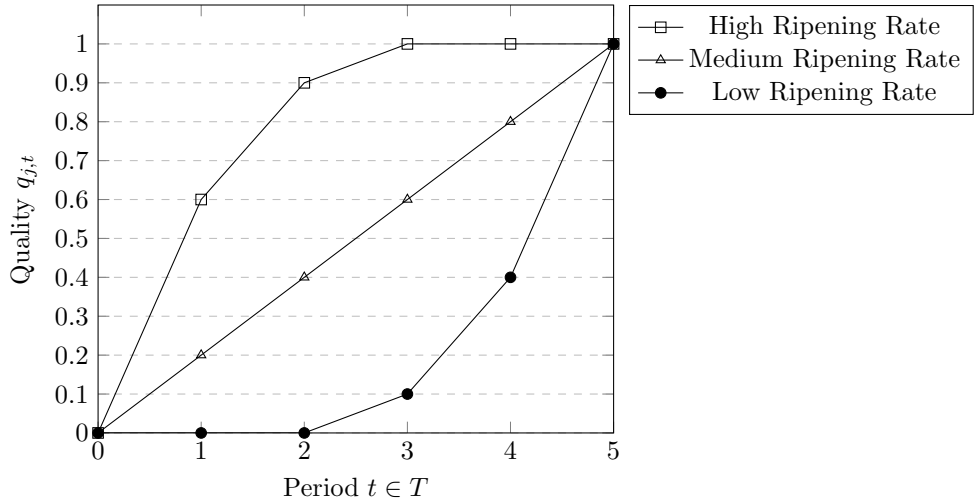

Then, for the second analysis to simulate having different ripening rates of quality we considered grape blocks whose qualities all start from zero and finish in one, but having different improvement rate. Graph \ref{q2} shows a representation of this case.

\section{Results}

In this section, we will present the main results. Starting by the effect that errors in grape yields estimation have on different evaluation metrics and then proceed with the effect of errors in the forecast of transition probabilities. Each analysis results will show first the outcomes from the different maximum reachable quality case and then, the ones from the different ripening rates of grapes case.

The selected metrics are: Objective function values, Objective function nominal and deficit percentage, Net income, the Proportion of total incomes increment/loss given by two different factors: The extra or less unharvested grape compared to the perfect information scenario and the average quality of harvested grape. These results will present the performance separately by the three levels of flexibility considered. It will be also showed the Percentage unharvested grape and the Loss income, distinctly for the different grape quality behaviors. As these results are similar for the three levels of flexibility except that the difference between the different grapes increases, we will present just the corresponding to flexibility level 2.

These results will let us compare, analyze and thus, identify the economic and production performance impacts of the harvest plans as results of the forecast accuracy.

%{\color{blue} Seria bueno hacer una intro un poco mas extensa de los resultados, para explicitar lo que se busca presentar}

\subsection{Grape Yields Analysis}

\subsubsection{Different Maximum Reachable Quality Case}
 
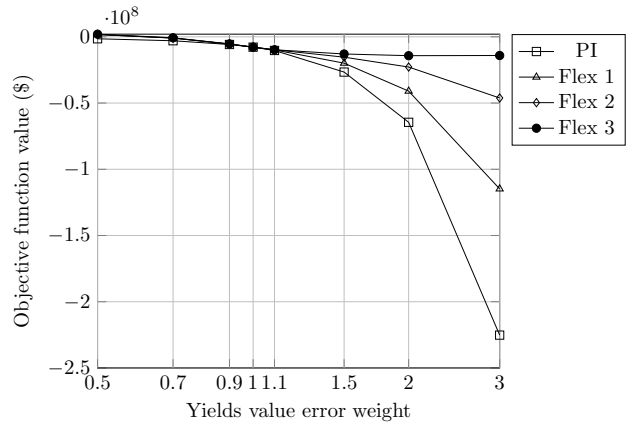
\begin{figure} [H]
\begin{center}
\begin{tikzpicture} [scale = 0.78]
\begin{semilogxaxis}[
    xlabel={Yields value error weight},
    ylabel={Objective function value (\$)},
    xmin=0.5, xmax=3,
    ymin=-250000000, ymax=2000000,
    xtick=\empty,
    extra x ticks      ={0.5,0.7,0.9,1,1.1,1.5,2,3},
    extra x tick labels={0.5,0.7,0.9,1,1.1,1.5,2,3},
    legend pos = outer north east,
    ymajorgrids=true,
    grid, ]

\addplot[mark=square] coordinates             
    {(3,-225199575)(2,-64539765)(1.5,-26615266)(1.1,-10311100)
    (1,-7809863)(0.9,-5793204)(0.7,-2974322)(0.5,-1487192)};

\addplot[mark=triangle] coordinates 
    {(3,-114811405)(2,-40992218)(1.5,-19887394)
    (1.1,-9885190)(1,-7809863)(0.9,-5457047)(0.7,-953122)(0.5,1704129)};
    
\addplot [mark=diamond] coordinates
    {(3,-46180260)(2,-22815230)(1.5,-15368651)(1.1,-9770914)
    (1,-7809863)(0.9,-5457047)(0.7,-945498)(0.5,1723111)};
    
\addplot[mark=*] coordinates 
    {(3,-14107332)(2,-14238771)(1.5,-12965651)
    (1.1,-9770914)(1,-7809863)(0.9,-5369304)(0.7,-726487)(0.5,1973337)};
    
\legend{PI, Flex 1, Flex 2,Flex 3},

\end{semilogxaxis}
\end{tikzpicture}
%%%%%%%%%%%%%%%%%%%%%%%%%%%%%%%%%%%%%%%%%%%%%%%%%%%%%%%%%%%%
%\begin{tikzpicture} [scale = 0.78]
%\begin{semilogxaxis}[
%    title={Objective function deficit},
%    xlabel={Yields value error weight},
%    ylabel={Objective function deficit (\$)},
%    xmin=0.5, xmax=3,
%    ymin=0, ymax=220000000,
%    xtick=\empty,
%    extra x ticks      ={0.5,0.7,0.9,1,1.1,1.5,2,3},
%    extra x tick labels={0.5,0.7,0.9,1,1.1,1.5,2,3},
%    legend pos = outer north east,
%    ymajorgrids=true,
%    grid, ]

%    \addplot[mark=square] coordinates 
%    {(3,110388170)(2,23547547)(1.5,6727872)(1.1,425910)(1,0)
%    (0.9,336157 )(0.7,2021200)(0.5,3191321)};
            
%    \addplot [mark=*] coordinates 
%    {(3,179019315)(2,41724535)(1.5,11246615)(1.1,425910)(1,0)(0.9,336157)(0.7,2028824)(0.5,3210303)};
        
%    \addplot [mark=triangle] coordinates 
%    {(3,211092243)(2,50300994)(1.5,13649615)(1.1,851739)(1,0)(0.9,423900)(0.7,2247835)(0.5,3460529)};
    
%    \legend{Flex 1, Flex 2,Flex 3},%

%\end{semilogxaxis}
%\end{tikzpicture}
\caption{Objective function differences in value for grape yields errors and flexibility levels}
\label{yield_objv}
\end{center}
\end{figure}

Figure \ref{yield_objv} show the change in the objective function value under different grape yields errors, for three levels of flexibility to modify the labor hiring schedule. We can observe that the effect of the errors in yield determination are not symmetrical; underestimations of the yields have a more significant, even negative effect in the objective function, while overestimation does not have such a significant effect. Flexibility can significantly affect the objective value, as is increased which means the decision-maker has the possibility to fix its decisions in an earlier period, the objective function change is smaller, so is less affected by errors in the beliefs. If yields are underestimated, flexibility plays an important role since the objective value is reduced exponentially as the error is increased and as the flexibility is increased, this effect is ameliorated.  
%The second graph shows specifically the objective function deficit compared to the perfect information scenario. Both graphs evidence that underestimating the grape yield and hence, the future stock, implicates high income losses while overestimating it does not.
%%%%%%%%%%%%%%%%%%%%%%%%%%%%%%%%%%%%%%%%%%%%%%%%%%%%%%%%%%%%

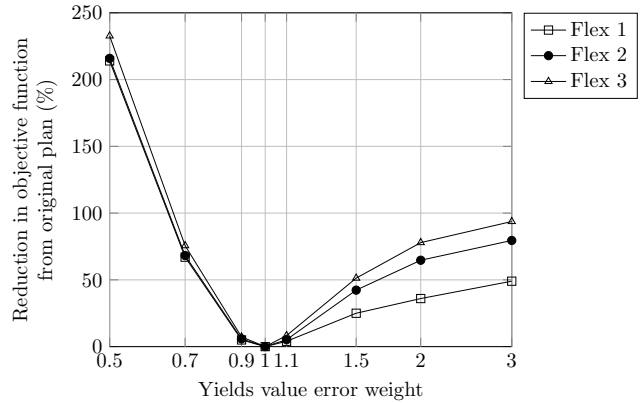
\begin{figure} [H]
\begin{center}
\begin{tikzpicture}[scale = 0.78]
\begin{semilogxaxis}[
    xlabel={Yields value error weight},
    ylabel style={align=center},
    ylabel =Reduction in objective function\\ from original plan (\%),
    xmin=0.5, xmax=3,
    ymin=0, ymax=250,
    xtick=\empty,
    extra x ticks      ={0.5,0.7,0.9,1,1.1,1.5,2,3},
    extra x tick labels={0.5,0.7,0.9,1,1.1,1.5,2,3},
    legend pos = outer north east,
    ymajorgrids=true,
    grid ]

    \addplot[mark=square] coordinates {(3,49)(2,36)(1.5,25)(1.1,4)(1,0)(0.9,5)(0.7,67)(0.5,214)};
            
    \addplot [mark=*] coordinates {(3,79.5)(2,64.7)(1.5,42.3)(1.1,5.2)(1,0)(0.9,5.8)(0.7,68.21)(0.5,215.9)};
        
    \addplot [mark=triangle] coordinates {(3,93.7)(2,77.9)(1.5,51.3)(1.1,8.3)(1,0)(0.9,7.3)(0.7,75.6)(0.5,232.7)};
    
    \legend{Flex 1, Flex 2,Flex 3},

\end{semilogxaxis}
\end{tikzpicture}
\caption{Percentage reduction in objective function from original plan (No error in the beliefs) by flexibility level.}
\label{yield_objp}
\end{center}
\end{figure}

\begin{figure} [H]
\begin{center}
%%%%%%%%%%%%%%%%%%%%%%%%%%%%%%%%%%%%%%%%%%%%%%%%%%%%%%%%%%%%
\begin{tikzpicture}[scale = 0.78]
\begin{semilogxaxis}[
    xlabel={Yields value error weight},
    ylabel style={align=center},
    ylabel=Reduction in income from\\ original plan (\%),
    xmin=0.5, xmax=3,
    ymin=-20, ymax=100,
    xtick=\empty,
    extra x ticks      ={0.5,0.7,0.9,1,1.1,1.5,2,3},
    extra x tick labels={0.5,0.7,0.9,1,1.1,1.5,2,3},
    legend pos = outer north east,
    ymajorgrids=true,
    grid ]

    \addplot[mark=square] coordinates {(3,40.78)(2,18.24)(1.5,24.97)(1.1,4.71)(1,0)(0.9,-5.24)(0.7,-10.55)(0.5,-17.41)};
            
    \addplot [mark=*] coordinates {(3,63.48)(2,46.38)(1.5,45.73)(1.1,9.05)(1,0)(0.9,-5.24)(0.7,-10.78)(0.5,-17.63)};
    
    \addplot [mark=triangle] coordinates {(3,95.96)(2,83.68)(1.5,62.49)(1.1,14.28)(1,0)(0.9,-6.32)(0.7,-10.99)(0.5,-17.72)};
    
    \legend{Flex 1, Flex 2,Flex 3},
\end{semilogxaxis}
\end{tikzpicture}
\caption{Percentage reduction in income from original plan (No error in the beliefs) by flexibility level.}
\label{yield_obji}
\end{center}
\end{figure}
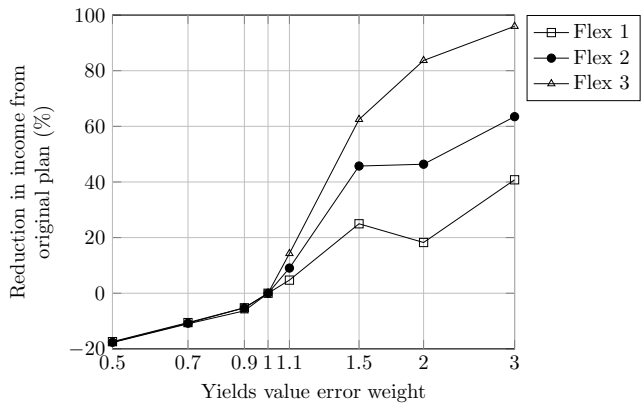

Figures \ref{yield_objp} and \ref{yield_obji} show the percentage reduction, with respect to the plan with the original belief (No error), in the objective value and income. The reduction in the objective function, as a percentage of the original plan, is much larger when the yields are overestimated with up to a 200\% loss in the value when the yields are 50\% of what it was expected. In the case when the yields are underestimated, when the underestimated level is in the 100\%, the loss in the objective value can range from 36\% to 77.9\%. If we observe the income, as expected, when the yields are overestimated they are reduced and when they are underestimated, they are increased compared to the original plan. However, neither the reduction or increase in income is proportional to the percentage reduction or increment. When yields are lowered by 50\% the incomes are only reduced by 17\%. In the case when yields are underestimated by 100\%, the income is incremented between 18\% to 83\%. 

Flexibility, defined as the ability to revise the hiring decisions, does not make a significant difference  in either the objective function or profits if the yields are overestimated. However, when the yields are underestimated flexibility plays an important role in reducing the effect over the objective value and the income. Having the ability to immediately re-plan in the second period can reduce the objective function and income by 36\% and 25\% respectively, when the yields are underestimated by 100\%. If flexibility is reduced by allowing only to re-plan in the fourth period, the objective function and income is now reduced by 77.9\% and 83.68\% respectively, when the yields are underestimated by 100\%.
%%%%%%%%%%%%%%%%%%%%%%%%%%%%%%%%%%%%%%%%%%%%%%%%%%%%%%%%%%%%

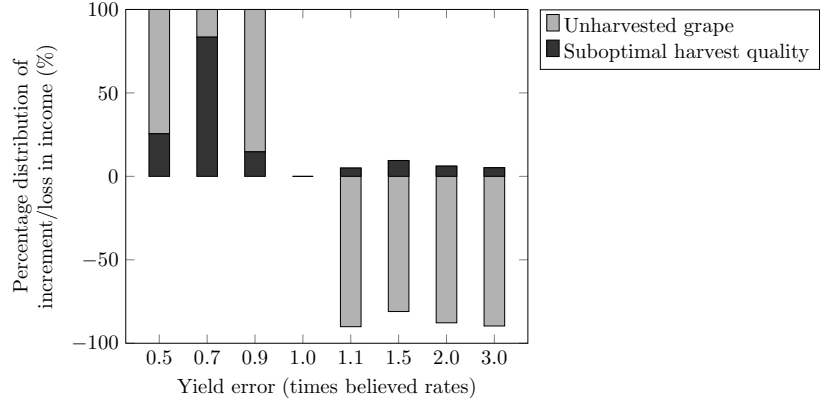
\begin{figure} [H]
\begin{center}
\begin{tikzpicture} [scale = 0.78]
\pgfplotsset{compat=1.9,compat/bar nodes=1.8,}
    \pgfplotstableread{
        Label series1 series2 
        0.5  25.53   74.47  
        0.7  83.40   16.60  
        0.9  14.69   85.31  
        1.0  00.00   00.00
        1.1   4.96  -95.04  
        1.5   9.50  -90.50  
        2.0   6.12  -93.88  
        3.0   5.17  -94.83  
    }\testdata
    \begin{axis}[
        xlabel={Yield error (times believed rates)},
        ylabel style={align=center},
        ylabel=Percentage distribution of\\ increment/loss in income (\%),
        ybar stacked,
        ymin=-100,
        ymax=100,
        xtick=data,
        legend style={
            cells={anchor=west},
            legend pos= outer north east,
        },
        reverse legend=true,
        xticklabels from table={\testdata}{Label},
        xticklabel style={text width=2cm,align=center},
    ]
        \addplot [fill=black!80]
            table [y=series1, meta=Label, x expr=\coordindex]
                {\testdata};
                    \addlegendentry{Suboptimal harvest quality}
        \addplot [fill=gray!60]
            table [y=series2, meta=Label, x expr=\coordindex]
                {\testdata};
                    \addlegendentry{Unharvested grape}
    \end{axis}
\end{tikzpicture}
\caption{Distribution of income increment/loss between unharvested and sub-optimally harvested grapes.}
\label{yield_objd}
\end{center}
\end{figure}

In Figure \ref{yield_objd} we can observe how the income increment/loss is distributed between unharvested and sub-optimally harvested grapes, compared to the case in which there was no error in the beliefs, for different levels of errors in the beliefs. We can observe that in the case when the yields are underestimated, the main source of income losses comes from the increasing levels of unharvested grape. In the other hand, overestimating these rates implicate an increment of harvested grape and also a higher average quality of it, which means more profits per unit. This compensates part of what is lost by the over-hiring labor.
%%%%%%%%%%%%%%%%%%%%%%%%%%%%%%%%%%%%%%%%%%%%%%%%%%%%%%%%%%%%%%%%%%%%%%
%%%%%%%%%%%%%%%%%%%%%%%%%%%%%%%%%%%%%%%%%%%%%%%%%%%%%%%%%%%%%%%%%%%%

\begin{figure} [H]
\begin{center}
\begin{tikzpicture}[scale = 0.78]
\begin{semilogxaxis}[
    xlabel={Yield error (times believed rates)},
    ylabel={Unharvested grape (\%)},
    xmin=0.5, xmax=3,
    ymin=-60, ymax=100,
    xtick=\empty,
    extra x ticks      ={0.5,0.7,0.9,1,1.1,1.5,2,3},
    extra x tick labels={0.5,0.7,0.9,1,1.1,1.5,2,3},
    legend pos = outer north east,
    ymajorgrids=true,
    grid, ]

    \addplot[mark=square] coordinates {(3,63.80)(2,45.31)(1.5,40.46)(1.1,5.49)(1,0)(0.9,1.36)(0.7,-1.49)(0.5,14.85)};
            
    \addplot [mark=*] coordinates {(3,69.00)(2,52.64)(1.5,58.46)(1.1,16.56)(1,0)(0.9,-6.78)(0.7,0.06)(0.5,-13.79)};
    
    \addplot [mark=triangle] coordinates {(3,91.28)(2,74.79)(1.5,67.95)(1.1,7.56)(1,0)(0.9,-10.36)(0.7,-2.46)(0.5,-57.38)};
    
    \legend{High Maximum Quality, Medium Maximum Quality, Low Maximum Quality},

\end{semilogxaxis}
\end{tikzpicture}
\caption{Percentage unharvested grape by maximum reachable quality.}
\label{yield_unharv1}
\end{center}
\end{figure}
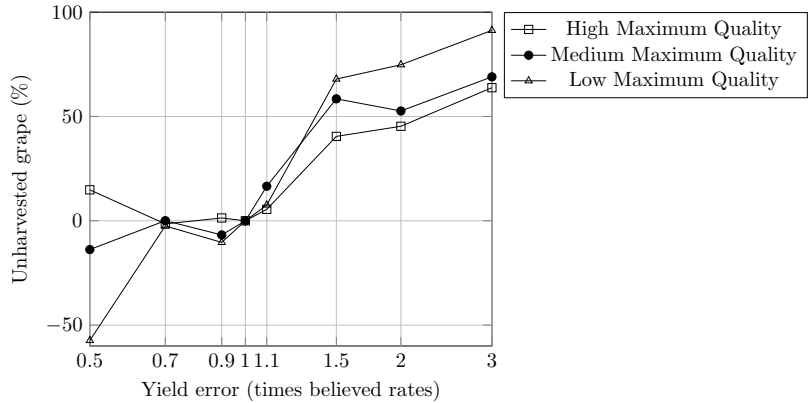

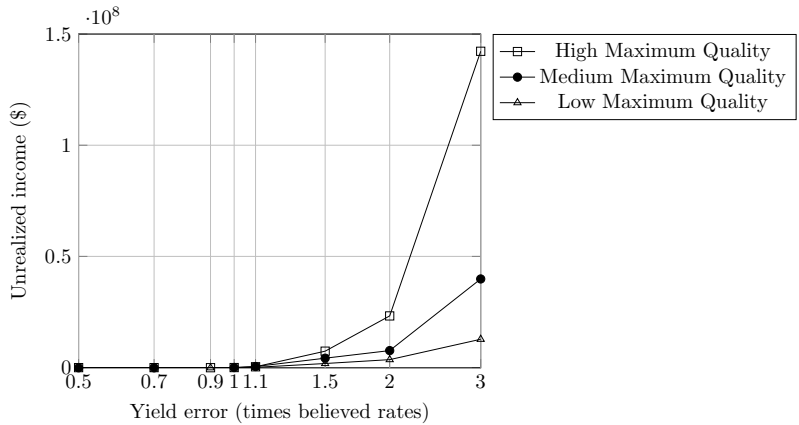
\begin{figure} [H]
\begin{center}
\begin{tikzpicture}[scale = 0.78]
\begin{semilogxaxis}[
    xlabel={Yield error (times believed rates)},
    ylabel={Unrealized income (\$)},
    xmin=0.5, xmax=3,
    ymin=-130000, ymax=150000000,
    xtick=\empty,
    extra x ticks      ={0.5,0.7,0.9,1,1.1,1.5,2,3},
    extra x tick labels={0.5,0.7,0.9,1,1.1,1.5,2,3},
    legend pos = outer north east,
    ymajorgrids=true,
    grid, ]

    \addplot[mark=square] coordinates {(3,142253456)(2,23231517)(1.5,7370949)(1.1,351790)(1,0)(0.9,-56825)(0.7,-24981)(0.5,-24266)};
            
    \addplot [mark=*] coordinates {(3,39850630)(2,7633685)(1.5,4215741)(1.1,472765)(1,0)(0.9,-143351)(0.7,-122103)(0.5,-91516)};
    
    \addplot [mark=triangle] coordinates {(3,12729371)(2,3565581)(1.5,1821464)(1.1,97641)(1,0)(0.9,-95992)(0.7,-163807)(0.5,-120996)};
    
    \legend{High Maximum Quality, Medium Maximum Quality, Low Maximum Quality},

\end{semilogxaxis}
\end{tikzpicture}
\caption{Unrealized income by maximum reachable quality.}
\label{yield_unreal1}
\end{center}
\end{figure}

Figures \ref{yield_unharv1} and \ref{yield_unreal1} present the results for three distinct blocks with different maximum reachable quality of grapes for flexibility two. We can observe that in the case of an underestimation scenario, the model privileges the most premium grape (high maximum quality) trying to leave unharvested the less amount as possible, while the worst grape suffers high percentage of loses in the face of variability. If we look at the effect on income,  we can notice that even though there is a lower level of losses of the better quality grapes, since they correspond to a significant proportion of the incomes, they account for the largest portion of the income loss. Hence the decision-maker must prioritize the premium blocks allocating labor to these blocks as much as the firm can. 

%%%%%%%%%%%%%%%%%%%%%%%%%%%%%%%%%%%%%%%%%%%%%%%%%%%%%%%%%%%%%%%%%%

\subsubsection{Different Ripening Rates of Quality Case}

%{\color{blue} Bajo qué condiciones estan los graficos de abajo}

\begin{figure} [H]
\begin{center}
\begin{tikzpicture}[scale = 0.78]
\begin{semilogxaxis}[
    xlabel={Yield error (times believed rates)},
    ylabel={Objective function value (\$)},
    xmin=0.5, xmax=3,
    ymin= -3000000000, ymax=15000000,
    xtick=\empty,
    extra x ticks      ={0.5,0.7,0.9,1,1.1,1.5,2,3},
    extra x tick labels={0.5,0.7,0.9,1,1.1,1.5,2,3},
    legend pos = outer north east,
    ymajorgrids=true,
    grid, ]

    \addplot[mark=square] coordinates {(3,-2978082055)(2,-486355485)(1.5,-134843415)(1.1,-34554162)(1,-22966637)(0.9,-14838339)(0.7,-5463596)(0.5,-1818274)};
            
    \addplot [mark=*] coordinates {(3,-1248080902)(2,-239295716)(1.5,-81115481)(1.1,-31452484)(1,-22966637)(0.9,-12170428)(0.7,7545825)(0.5,14327191)};
    
    \addplot [mark=triangle] coordinates {(3,-274071491)(2,-90461863)(1.5,-51766016)(1.1,-31106091)(1,-22966637)(0.9,-12169396)(0.7,7705111)(0.5,14636090)};
    
    \addplot [mark=diamond] coordinates {(3,-80417796)(2,-66871955)(1.5,-51310809)(1.1,-31106091)(1,-22966637)(0.9,-12170428)(0.7,7863799)(0.5,14950233)};
    
    \legend{PI , Flex 1, Flex 2, Flex 3},

\end{semilogxaxis}
\end{tikzpicture}
\caption{Objective function differences in value for grape yields errors and flexibility levels.}
\label{rates_obj1.2}
\end{center}
\end{figure}
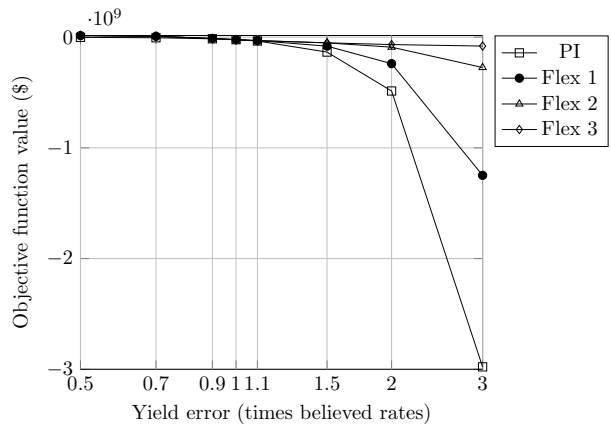

In the figure \ref{rates_obj1.2} we present the aggregate performance of harvest plans which deals with grapes with different ripening rates under the three different levels of flexibility considered. Here can observe that the mentioned flexibility to fix the harvest labor plan has similar impacts on the objective function. Having the possibility to fix the decisions immediately after, in the second period, manages to reduce the losses in about 40\%. 

%{\color{blue} Ver las leyendas de estos gráficos, ya que no indican cual es la tasa de maduración, sale felxibilidad. O esto es para una tasa dada... no queda claro que es lo que estamos viendo}

%%%%%%%%%%%%%%%%%%%%%%%%%%%%%%%%%%%%%%%%%%%%%%%%%%%%%%%%%%%%%%%%%%%%
%%%%%%%%%%%%%%%%%%%%%%%%%%%%%%%%%%%%%%%%%%%%%%%%%%%%%%%%%%%%%%%%%%%%%%%%%%%

%\begin{tikzpicture}[scale = 0.78]
%\begin{semilogxaxis}[
%    title={Objective function deficit},
%    xlabel={Yield error (times believed rates)},
%    ylabel={Objective function deficit (\$)},
%    xmin=0.5, xmax=3,
%    ymin=0, ymax=3000000000,
%    xtick=\empty,
%    extra x ticks      ={0.5,0.7,0.9,1,1.1,1.5,2,3},
%    extra x tick labels={0.5,0.7,0.9,1,1.1,1.5,2,3},
%    legend pos = outer north east,
%    ymajorgrids=true,
%    grid, ]

%    \addplot[mark=square] coordinates 
%    {(3,1730001153)(2,247059769)(1.5,53727934)(1.1,3101678)(1,0)(0.9,2667911)(0.7,13009421)(0.5,16145465)};
            
%    \addplot [mark=*] coordinates 
%    {(3,2704010564)(2,395893622)(1.5,83077399)(1.1,3448071)(1,0)(0.9,2668943)(0.7,13168707)(0.5,16454364)};
        
%    \addplot [mark=triangle] coordinates 
%    {(3,2897664259)(2,419483530)(1.5,83532606)(1.1,3448071)(1,0)(0.9,2667911)(0.7,13327395)(0.5,16768507)};
    
%    \legend{Flex 1, Flex 2,Flex 3},

%\end{semilogxaxis}
%\end{tikzpicture}

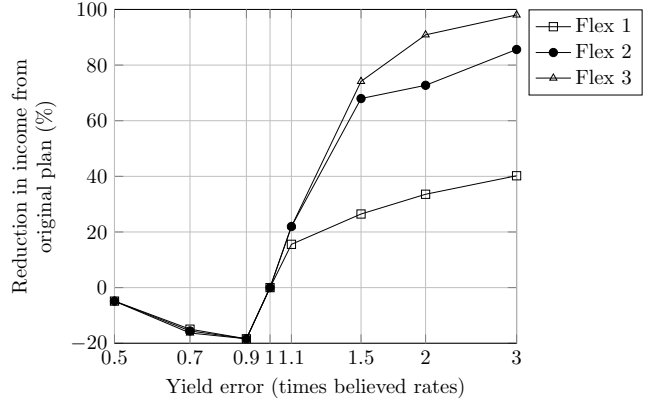
\begin{figure} [H]
\begin{center}
\begin{tikzpicture}[scale = 0.78]
\begin{semilogxaxis}[
    xlabel={Yield error (times believed rates)},
    ylabel style={align=center},
    ylabel=Reduction in income from\\ original plan (\%),
    xmin=0.5, xmax=3,
    ymin=-20, ymax=100,
    xtick=\empty,
    extra x ticks      ={0.5,0.7,0.9,1,1.1,1.5,2,3},
    extra x tick labels={0.5,0.7,0.9,1,1.1,1.5,2,3},
    legend pos = outer north east,
    ymajorgrids=true,
    grid, ]

    \addplot[mark=square] coordinates {(3,40.21)(2,33.56)(1.5,26.46)(1.1,15.56)(1,0)(0.9,-18.4)(0.7,-14.94)(0.5,-4.85)};
            
    \addplot [mark=*] coordinates {(3,85.61)(2,72.73)(1.5,67.97)(1.1,21.96)(1,0)(0.9,-18.4)(0.7,-15.57)(0.5,-4.85)};
    
    \addplot [mark=triangle] coordinates {(3,98.02)(2,90.87)(1.5,74.13)(1.1,21.96)(1,0)(0.9,-18.4)(0.7,-16.31)(0.5,-4.86)};
    
    \legend{Flex 1, Flex 2,Flex 3},

\end{semilogxaxis}
\end{tikzpicture}
\caption{Percentage reduction in income from original plan (No error in the beliefs) for different flexibility levels.}
\label{rates_income1.2}
\end{center}
\end{figure}

We can observe in figure \ref{rates_income1.2}, similar impacts on net incomes the fact of under and over estimating grape yields with different ripening rates of grape. About 60\% of net income losses can be avoided if investment on improving flexibility are made. In the next figure, we can notice that additional unharvested grape is the main source of income losses, while when rates are overestimated, the extra incomes sources are more equally generated by both factors. 
%%%%%%%%%%%%%%%%%%%%%%%%%%%%%%%%%%%%%%%%%%%%%%%%%%%%%%%%%%%%%%%%%%%%%%%%%%%

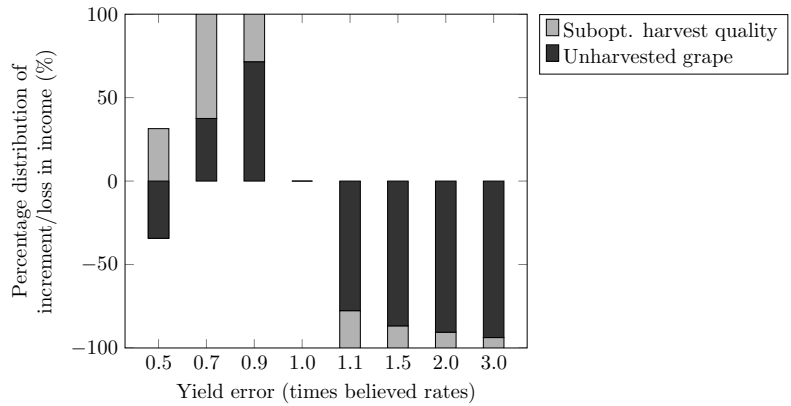
\begin{figure} [H]
\begin{center}
%%%%%%%%%%%%%%%%%%%%%%%%%%%%%%%%%%%%%%%%%%%%%%%%%%%%%%%%%%%%%%%%%%%%%%
\begin{tikzpicture}[scale = 0.78]
\pgfplotsset{compat=1.9,compat/bar nodes=1.8,}
    \pgfplotstableread{
        Label series1 series2 
        0.5  -34.29   65.71  
        0.7   37.50   62.50  
        0.9   71.40   28.60  
        1.0   00.00   00.00
        1.1  -77.66  -22.34  
        1.5  -86.75  -13.25  
        2.0  -90.57  -09.43   
        3.0  -93.77  -06.23   
    }\testdata
    \begin{axis}[
        xlabel={Yield error (times believed rates)},
        ylabel style={align=center},
        ylabel=Percentage distribution of\\ increment/loss in income (\%),
        ybar stacked,
        ymin=-100,
        ymax=100,
        xtick=data,
        legend style={
            cells={anchor=west},
            legend pos= outer north east,
        },
        reverse legend=true,
        xticklabels from table={\testdata}{Label},
        xticklabel style={text width=2cm,align=center},
    ]
        \addplot [fill=black!80]
            table [y=series1, meta=Label, x expr=\coordindex]
                {\testdata};
                    \addlegendentry{Unharvested grape}
        \addplot [fill=gray!60]
            table [y=series2, meta=Label, x expr=\coordindex]
                {\testdata};
                    \addlegendentry{Subopt. harvest quality}
    \end{axis}
\end{tikzpicture}
\caption{Distribution of income increment/loss between unharvested and sub-optimally harvested grapes.}
\label{rates_distrib}
\end{center}
\end{figure}

In the following graphs we can observe a very important behavior of each quality behavior. The results show that grapes that have an early improvement of their quality gives the decision-maker an extra level of flexibility to adjust the harvesting plan. As this kind of grape reaches its optimal quality early and stay in that condition for more periods, the planner is able to start its harvest if necessary. This is evidenced in the lower amount of unharvested grape as the error increases. In the other hand, we can see that grapes which have a late quality improvement limits the flexibility of the decision-maker to extend the harvest period due to the quality does not make it profitable. These results are almost the same for the other two levels of flexibility an so, we presented just the corresponding to flexibility level 2.
%%%%%%%%%%%%%%%%%%%%%%%%%%%%%%%%%%%%%%%%%%%%%%%%%%%%%%%%%%%%%%%%%%%%%%

%{\color{blue} Para qué nivel de flexibilidad esta definido el grafico de abajo?}

\begin{figure} [H]
\begin{center}
\begin{tikzpicture}[scale = 0.78]
\begin{semilogxaxis}[
    xlabel={Yield error (times believed rates)},
    ylabel={Unharvested grape (\%)},
    xmin=0.5, xmax=3,
    ymin=-22, ymax=100,
    xtick=\empty,
    extra x ticks      ={0.5,0.7,0.9,1,1.1,1.5,2,3},
    extra x tick labels={0.5,0.7,0.9,1,1.1,1.5,2,3},
    legend pos = outer north east,
    ymajorgrids=true,
    grid, ]

    \addplot[mark=square] coordinates {(3,55.76)(2,35.97)(1.5,40.67)(1.1,20.97)(1,0)(0.9,-21.75)(0.7,-7.35)(0.5,6.32)};
            
    \addplot [mark=*] coordinates {(3,87.46)(2,72.05)(1.5,64.03)(1.1,13.59)(1,0)(0.9,-15.48)(0.7,-4.74)(0.5,3.91)};
    
    \addplot [mark=triangle] coordinates {(3,96.23)(2,84.45)(1.5,78.38)(1.1,17.65)(1,0)(0.9,-4.76)(0.7,-4.55)(0.5,-4.35)};
    
    \legend{High Ripening Rate, Medium Ripening Rate, Low Ripening Rate},

\end{semilogxaxis}
\end{tikzpicture}
\caption{Percentage unharvested grape for different ripening rates of grapes.}
\label{yield_unreal}
\end{center}
\end{figure}
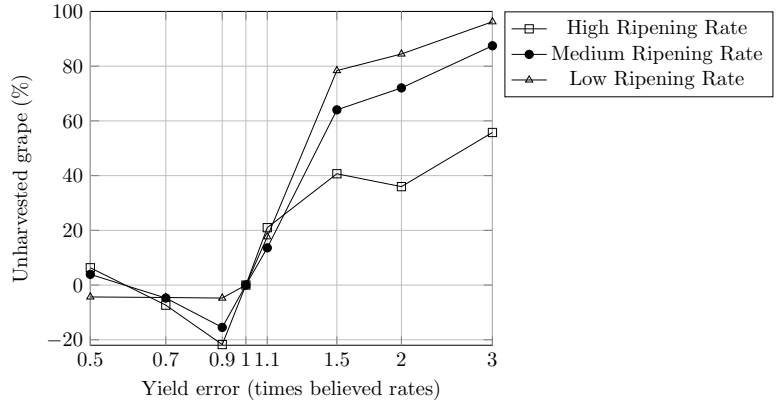

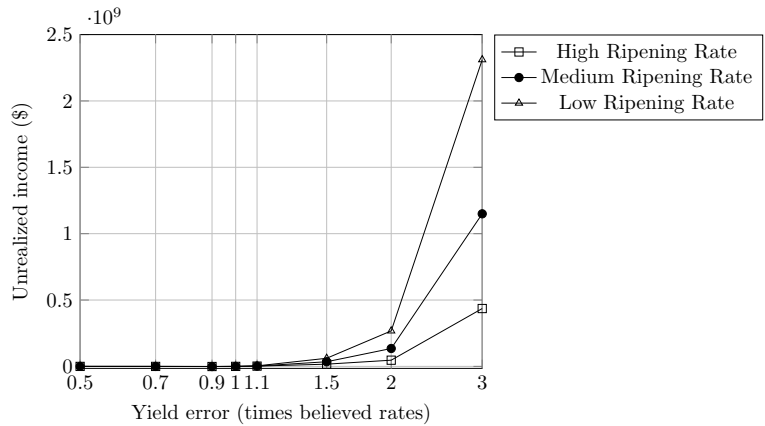
\begin{figure} [H]
\begin{center}
%%%%%%%%%%%%%%%%%%%%%%%%%%%%%%%%%%%%%%%%%%%%%%%%%%%%%%%%%%%%%%%%%%%%%%
\begin{tikzpicture}[scale = 0.78]
\begin{semilogxaxis}[
    xlabel={Yield error (times believed rates)},
    ylabel={Unrealized income (\$)},
    xmin=0.5, xmax=3,
    ymin=-15000000, ymax=2500000000,
    xtick=\empty,
    extra x ticks      ={0.5,0.7,0.9,1,1.1,1.5,2,3},
    extra x tick labels={0.5,0.7,0.9,1,1.1,1.5,2,3},
    legend pos = outer north east,
    ymajorgrids=true,
    grid, ]

    \addplot[mark=square] coordinates {(3,435881326)(2,46710509)(1.5,16788523)(1.1,2033993)(1,0)(0.9,-904258)(0.7,-169188)(0.5,-39770)};
            
    \addplot [mark=*] coordinates {(3,1149612062)(2,134820137)(1.5,35033950)(1.1,2193648)(1,0)(0.9,-1196656)(0.7,-419352)(0.5,-4505)};
    
    \addplot [mark=triangle] coordinates {(3,2310574505)(2,267379869)(1.5,61163273)(1.1,3918577)(1,0)(0.9,-551680)(0.7,-156418)(0.5,-23668)};
    
    \legend{High Ripening Rate, Medium Ripening Rate, Low Ripening Rate},

\end{semilogxaxis}
\end{tikzpicture}
\caption{Unrealized income for different ripening rates of grapes.}
\label{yield_unreal}
\end{center}
\end{figure}

\subsection{Transition Probabilities Analysis}

\subsubsection{Different Maximum Reachable Quality Case}
In this analysis we can observe a more linear behavior in comparison with the grape yield analysis. Despite this, the following graphs still evidence an asymmetric performance between overestimating the probabilities of good scenarios and underestimating them. Underestimating the probabilities of good scenarios can cost five times in comparison to overestimate them.
%%%%%%%%%%%%%%%%%%%%%%%%%%%%%%%%%%%%%%%%%%%%%%%%%%%%%%%%%%%%%%%%%%%%%
\begin{figure} [H]
\begin{center}
\begin{tikzpicture}[scale = 0.78]
\begin{axis}[
    ylabel={Objective function value (\$)},
    ymin=-25000000, ymax=0,
    xlabel={Scenarios},
    xmin=1, xmax= 8,
    xtick=\empty,
    extra x ticks      ={1,2,3,4,5,6,7,8},
    extra x tick labels={1,2,3,4,5,6,7,8},
    legend pos = outer north east,
    ymajorgrids=true,
    grid, ]

\addplot[mark=square] coordinates             
    {(8,-23570000)(7,-13997942)(6,-11204776)(5,-8263263)(4,-7809863)(3,-6094902)(2,-5671231)(1,-3834112)};

\addplot[mark=triangle] coordinates 
    {(8,-19859200)(7,-12628080)(6,-10446579)(5,-8181816)(4,-7809863)(3,-6027361)(2,-5067475)(1,-1858263)};
    
\addplot [mark=diamond] coordinates
    {(8,-16298400)(7,-12218671)(6,-10339694)(5,-8181665)(4,-7809863)(3,-6027413)(2,-5067475)(1,-1856835)};
    
\addplot[mark=*] coordinates 
    {(8,-16174300)(7,-12089024)(6,-10170420)(5,-8094791)(4,-7809863)(3,-6027049)(2,-5067475)(1,-1856835)};
    
\legend{PI, Flex 1, Flex 2,Flex 3},

\end{axis}
\end{tikzpicture}
\caption{Objective function value for different scenarios and flexibility levels.}
\label{yield_unreal}
\end{center}
\end{figure}
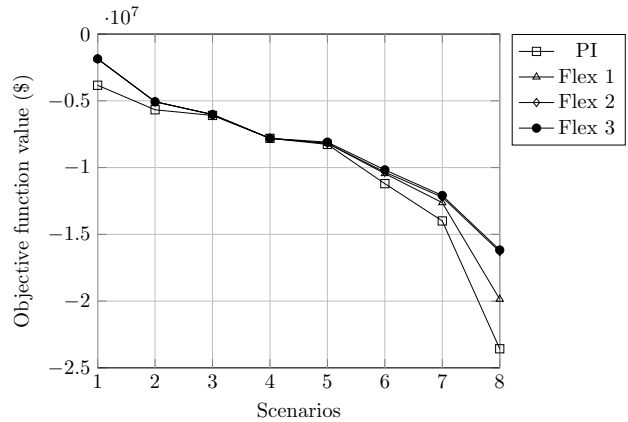

%%%%%%%%%%%%%%%%%%%%%%%%%%%%%%%%%%%%%%%%%%%%%%%%%%%%%%%%%%%%%%%%%%%%%%%
%\begin{tikzpicture}[scale = 0.78]
%\begin{axis}[
%    title={Objective function deficit},
%    ylabel={Objective function deficit (\$)},
%    ymin=0, ymax=52000000,
%    xlabel={1°stage decision error (\%)},
%    xmin=-63, xmax= 150,
%    xtick=\empty,
%    extra x ticks      ={-63,-33,-17,0,17,58,79,150},
%    extra x tick labels={-63,-33,-17,0,17,58,79,150},
%    legend pos = outer north east,
%    ymajorgrids=true,
%    grid, ]
    
%    \addplot[mark=square] coordinates 
%    {(150,21367500)(79,10378196)(58,4342044)(17,432253)(0,0)(-17,1095991)(-33,4489330)(-63,12848444)};
            
%    \addplot [mark=*] coordinates 
%    {(150,50242750)(79,13422146)(58,4939447)(17,432253)(0,0)(-17,1095698)(-33,4489330)(-63,13033944)};
        
%    \addplot [mark=triangle] coordinates 
%    {(150,50242750)(79,13422146)(58,4939534)(17,432253)(0,0)(-17,1095698)(-33,4489330)(-63,13219444)};
    
%    \legend{Flex 1, Flex 2,Flex 3},

%\end{axis}
%\end{tikzpicture}

In the graphs \ref{reduct2.1} and \ref{unharv2.1} we can notice that flexibility has less impact on global results than in the previous analysis. The reduction of deficit could reach just about a 15\%. In addition, we can observe that the less flexibility you have, the less impact has the improvement of flexibility to be able to fix the schedule one period before. A very important result that we can observe here, is the fact that flexibility makes no difference when bad scenarios probabilities are underestimated (left side of the graphs).
%%%%%%%%%%%%%%%%%%%%%%%%%%%%%%%%%%%%%%%%%%%%%%%%%%%%%%%%%%%%%%%%%%%%
\begin{figure} [H]
\begin{center}
\begin{tikzpicture}[scale = 0.78]
\begin{axis}[
    ylabel style={align=center},
    ylabel =Reduction in objective function\\ from original plan (\%),
    ymin=0, ymax=55,
    xlabel={Scenario},
    xmin=1, xmax= 8,
    xtick=\empty,
    extra x ticks      ={1,2,3,4,5,6,7,8},
    extra x tick labels={1,2,3,4,5,6,7,8},
    legend pos = outer north east,
    ymajorgrids=true,
    grid, ]

    \addplot[mark=square] coordinates {(8,15.74)(7,9.79)(6,6.77)(5,0.99)(4,0)(3,1.11)(2,10.65)(1,51.53)};
            
    \addplot [mark=*] coordinates 
    {(8,30.85)(7,12.71)(6,7.72)(5,0.99)(4,0)(3,1.11)(2,10.65)(1,51.57)};
        
    \addplot [mark=triangle] coordinates {(8,41.38)(7,13.64)(6,9.23)(5,2.04)(4,0)(3,1.11)(2,10.65)(1,51.57)};
    
    \legend{Flex 1, Flex 2,Flex 3},

\end{axis}
\end{tikzpicture}
\caption{Percentage reduction in objective function from original plan (No error in the beliefs) for different scenarios and flexibility levels.}
\label{reduct2.1}
\end{center}
\end{figure}
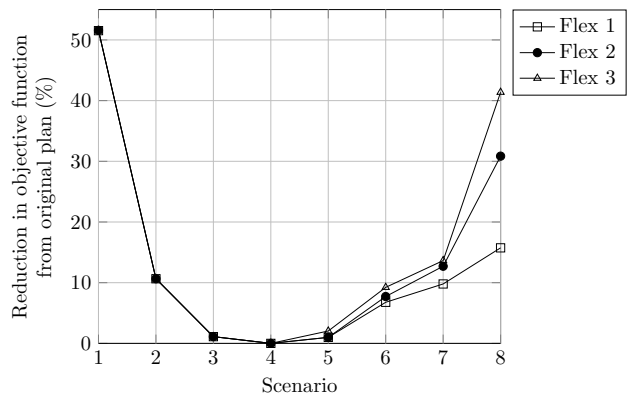

\begin{figure} [H]
\begin{center}
%%%%%%%%%%%%%%%%%%%%%%%%%%%%%%%%%%%%%%%%%%%%%%%%%%%%%%%%%%%%
\begin{tikzpicture}[scale = 0.78]
\begin{axis}[
    ylabel style={align=center},
    ylabel=Reduction in income from\\ original plan (\%),
    ymin=-13, ymax=35,
    xlabel={Scenario},
    xmin=1, xmax= 8,
    xtick=\empty,
    extra x ticks      ={1,2,3,4,5,6,7,8},
    extra x tick labels={1,2,3,4,5,6,7,8},
    legend pos = outer north east,
    ymajorgrids=true,
    grid, ]

    \addplot[mark=square] coordinates {(8,10.88)(7,10.56)(6,6.86)(5,4.88)(4,0)(3,-3.31)(2,-9.87)(1,-12.18)};
            
    \addplot [mark=*] coordinates {(8,28.08)(7,17.62)(6,13.01)(5,4.89)(4,0)(3,-3.37)(2,-9.87)(1,-12.37)};
    
    \addplot [mark=triangle] coordinates {(8,34.12)(7,19.98)(6,15.28)(5,5.02)(4,0)(3,-3.33)(2,-9.87)(1,-12.37)};
    
    \legend{Flex 1, Flex 2,Flex 3},

\end{axis}
\end{tikzpicture}
\caption{Percentage reduction in income from original plan (No error in the beliefs) for different scenarios and flexibility levels.}
\label{unharv2.1}
\end{center}
\end{figure}
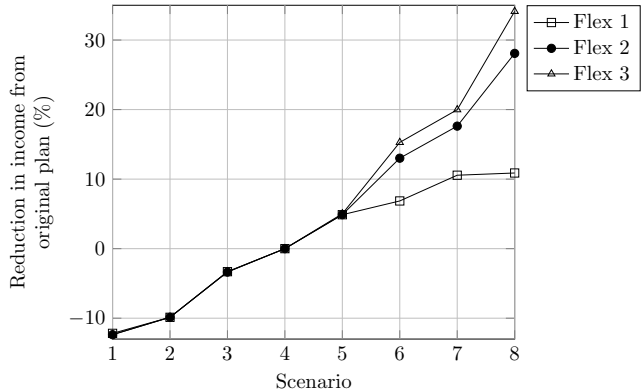

The next histogram shows, as in the previous analysis, the percentage implication on total incomes by two different factors: The extra or less unharvested grape compared to the perfect information scenario that means a variation on total incomes and penalty for leaving this amount unharvested, and the average quality of harvested grape that impacts on the revenues. We can observe similar results when yields are underestimated, the main source of income losses come from the increasing levels of unharvested grape too. But, overestimating these rates does not implicate the same as before, the main source of loss compensation is the increment of harvested grape. These results are almost the same for the other two levels of flexibility an so, we presented just the corresponding to flexibility level 2.
%%%%%%%%%%%%%%%%%%%%%%%%%%%%%%%%%%%%%%%%%%%%%%%%%%%%%%%%%%%%%%%%%%%%%%

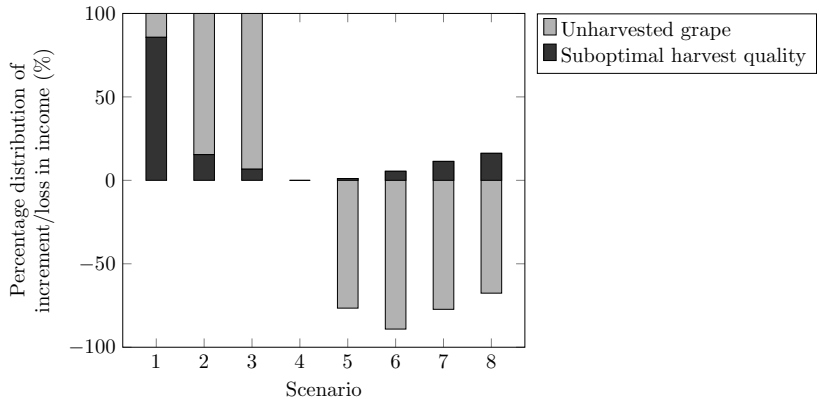
\begin{figure} [H]
\begin{center}
\begin{tikzpicture}[scale = 0.78]
\pgfplotsset{compat=1.9,compat/bar nodes=1.8,}
    \pgfplotstableread{
        Label series1 series2 
        1  85.70   14.30   
        2  15.41   84.59   
        3  06.68   93.32   
        4   00.00   00.00
        5   01.09  -77.66   
        6   05.44  -94.56   
        7   11.36  -88.64    
        8  16.20  -83.80    
    }\testdata
    \begin{axis}[
        xlabel={Scenario},
        ylabel style={align=center},
        ylabel=Percentage distribution of\\ increment/loss in income (\%),
        ybar stacked,
        ymin=-100,
        ymax=100,
        xtick=data,
        legend style={
            cells={anchor=west},
            legend pos= outer north east,
        },
        reverse legend=true,
        xticklabels from table={\testdata}{Label},
        xticklabel style={text width=2cm,align=center},
    ]
        \addplot [fill=black!80]
            table [y=series1, meta=Label, x expr=\coordindex]
                {\testdata};
                    \addlegendentry{Suboptimal harvest quality}
        \addplot [fill=gray!60]
            table [y=series2, meta=Label, x expr=\coordindex]
                {\testdata};
                    \addlegendentry{Unharvested grape}
    \end{axis}
\end{tikzpicture}
\caption{Distribution of income increment/loss between unharvested and sub-optimally harvested grapes for different scenarios.}
\label{yield_unreal}
\end{center}
\end{figure}

\begin{figure} [H]
\begin{center}
%%%%%%%%%%%%%%%%%%%%%%%%%%%%%%%%%%%%%%%%%%%%%%%%%%%%%%%%%%%%%%%%%%
\begin{tikzpicture}[scale = 0.78]
\begin{axis}[
    ylabel={Unrealized income (\$)},
    ymin=-250000, ymax=2700000,
    xlabel={Scenario},
    xmin=1, xmax= 8,
    xtick=\empty,
    extra x ticks      ={1,2,3,4,5,6,7,8},
    extra x tick labels={1,2,3,4,5,6,7,8},
    legend pos = outer north east,
    ymajorgrids=true,
    grid, ]

    \addplot[mark=square] coordinates {(8,2691800)(7,995495)(6,544130)(5,105088)(4,0)(3,-34752)(2,-63342)(1,-41378)};
            
    \addplot [mark=*] coordinates {(8,1874800)(7,827526)(6,670584)(5,247218)(4,0)(3,-111473)(2,-211778)(1,-188935)};
    
    \addplot [mark=triangle] coordinates {(8,1204000)(7,490290)(6,204672)(5,55900)(4,0)(3,-64653)(2,-237139)(1,-182876)};
    
    \legend{High Maximum Quality, Medium Maximum Quality, Low Maximum Quality},

\end{axis}
\end{tikzpicture}
\caption{Unrealized income for different scenarios and maximum reachable quality.}
\label{unreal2.1}
\end{center}
\end{figure}
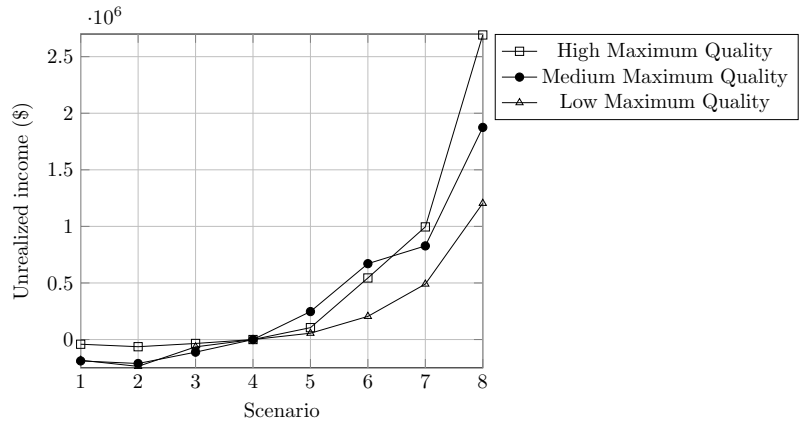

In the graphs \ref{unharv2.1} and \ref{unreal2.1} we present the main results separately by grapes with different maximum reachable quality for flexibility two, affected by the level of beliefs pessimism. In these, we observe similar outcomes. The most premium grape is prioritized as opposed to the lower quality grapes. Even so, it also happens that the extra unharvested premium quality grape composes the mayor part of the income losses.
%%%%%%%%%%%%%%%%%%%%%%%%%%%%%%%%%%%%%%%%%%%%%%%%%%%%%%%%%%%%%%%%%%%%%%
\begin{figure} [H]
\begin{center}
\begin{tikzpicture}[scale = 0.78]
\begin{axis}[
    ylabel={Unharvested grape (\%)},
    ymin=-20, ymax=35,
    xlabel={Scenario},
    xmin=1, xmax= 8,
    xtick=\empty,
    extra x ticks      ={1,2,3,4,5,6,7,8},
    extra x tick labels={1,2,3,4,5,6,7,8},
    legend pos = outer north east,
    ymajorgrids=true,
    grid, ]

    \addplot[mark=square] coordinates {(8,6.5)(7,6.08)(6,4.36)(5,2.63)(4,0)(3,-0.78)(2,3.03)(1,-2.1)};
            
    \addplot [mark=*] coordinates {(8,21.67)(7,15.87)(6,8.84)(5,9.01)(4,0)(3,-4.62)(2,-11.05)(1,-1.7)};
    
    \addplot [mark=triangle] coordinates {(8,32.5)(7,30.05)(6,14.99)(5,4.54)(4,0)(3,-4.99)(2,-19.78)(1,-5.69)};
    
    \legend{High Maximum Quality, Medium Maximum Quality, Low Maximum Quality},

\end{axis}
\end{tikzpicture}
\caption{Percentage unharvested grape for different scenarios and maximum reachable quality.}
\label{unharv2.1}
\end{center}
\end{figure}
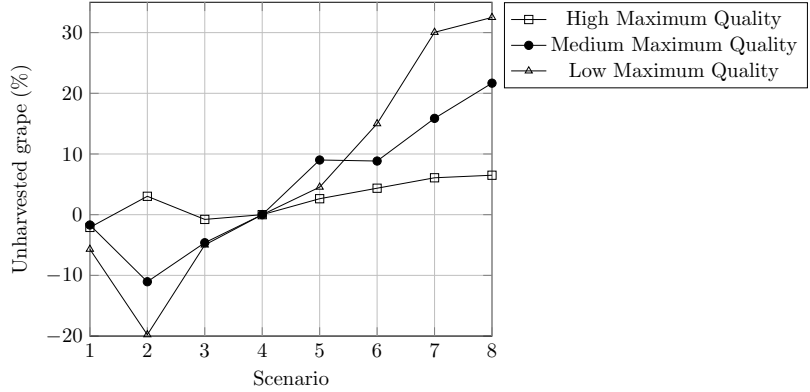

\subsubsection{Different ripening rates of grapes}

\begin{figure} [H]
\begin{center}
\begin{tikzpicture}[scale = 0.78]
\begin{axis}[
    ylabel={Objective function value (\$)},
    ymin=-120000000, ymax=7000000,
    xlabel={Scenario},
    xmin=1, xmax= 8,
    xtick=\empty,
    extra x ticks      ={1,2,3,4,5,6,7,8},
    extra x tick labels={1,2,3,4,5,6,7,8},
    legend pos = outer north east,
    ymajorgrids=true,
    grid, ]

\addplot[mark=square] coordinates             
    {(8,-109592000)(7,-54069641)(6,-38666944)(5,-25591245)
    (4,-22966637)(3,-15777766)(2,-14689517)(1,-8365594)};

\addplot[mark=triangle] coordinates 
    {(8,-88224500)(7,-43691445)(6,-34324900)(5,-25158992)
    (4,-22966637)(3,-14681775)(2,-10200187)(1,4482850)};
    
\addplot [mark=diamond] coordinates
    {(8,-69349250)(7,-40647495)(6,-33727497)(5,-25158992)
    (4,-22966637)(3,-14682068)(2,-10200187)(1,4668350)};
    
\addplot[mark=*] coordinates 
    {(8,-59349250)(7,-34647495)(6,-33727410)(5,-25158992)
    (4,-22966637)(3,-14682068)(2,-10200187)(1,4853850)};
    
\legend{PI, Flex 1, Flex 2,Flex 3},

\end{axis}
\end{tikzpicture}
\caption{Objective function value for different scenarios and flexibility levels.}
\label{yield_unreal2.2}
\end{center}
\end{figure}
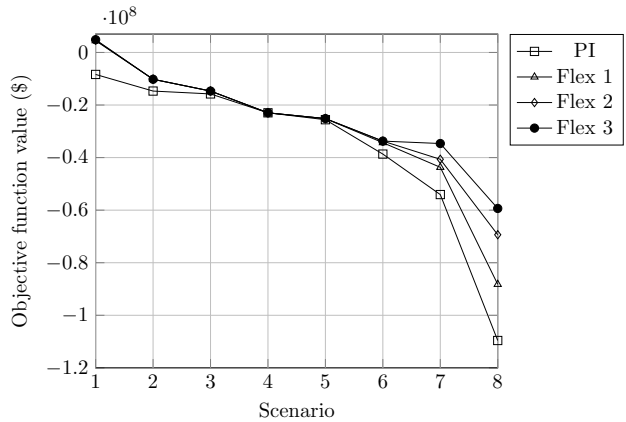

In the following graph (\ref{yield_unreal2.2}) we can observe that, in this case, flexibility makes no difference when bad scenarios have more occurrence probability than expected. On the other hand, flexibility plays again an important role when good scenarios probabilities are sub-estimated. 
%%%%%%%%%%%%%%%%%%%%%%%%%%%%%%%%%%%%%%%%%%%%%%%%%%%%%%%%%%%%%%%%%%%%
%%%%%%%%%%%%%%%%%%%%%%%%%%%%%%%%%%%%%%%%%%%%%%%%%%%%%%%%%%%%

%%%%%%%%%%%%%%%%%%%%%%%%%%%%%%%%%%%%%%%%%%%%%%%%%%%%%%%%%%%%%%%%%%%%
%\begin{figure} [H]
%\begin{center}
%\begin{tikzpicture}[scale = 0.78]
%\begin{axis}[
%    ylabel={Objective function value (\$)},
%    ymin=-80000000, ymax=3000000,
%    xlabel={Scenario},
%    xmin=1, xmax= 8,
%    xtick=\empty,
%    extra x ticks      ={1,2,3,4,5,6,7,8},
%    extra x tick labels={1,2,3,4,5,6,7,8},
%    legend pos = outer north east,
%    ymajorgrids=true,
%    grid, ]

%\addplot[mark=square] coordinates             
%    {(8,-75457700)(7,-48021294)(6,-39482253)(5,-29067591)
%    (4,-26474730)(3,-17788488)(2,-12717745)(1,1772350)};

%\addplot[mark=triangle] coordinates 
%    {(8,-50689700)(7,-33241273)(6,-27694604)(5,-21144981)
%    (4,-19989727)(3,-13436911)(2,-11112530)(1,1778565)};
    
%\addplot [mark=diamond] coordinates
%    {(8,-17094450)(7,-12502158)(6,-10756561)(5,-8897463)
%    (4,-9487810)(3,-6726651)(2,-7912625)(1,1778565)};
    
%\legend{High Ripening Rate, Medium Ripening Rate,Low Ripening Rate},

%\end{axis}
%\end{tikzpicture}
%\caption{Objective function value for different scenarios and ripening rates.}
%\label{yield_ripen2.2}
%\end{center}
%\end{figure}
%%%%%%%%%%%%%%%%%%%%%%%%%%%%%%%%%%%%%%%%%%%%%%%%%%%%%%%%%%%%%%%%%%%%
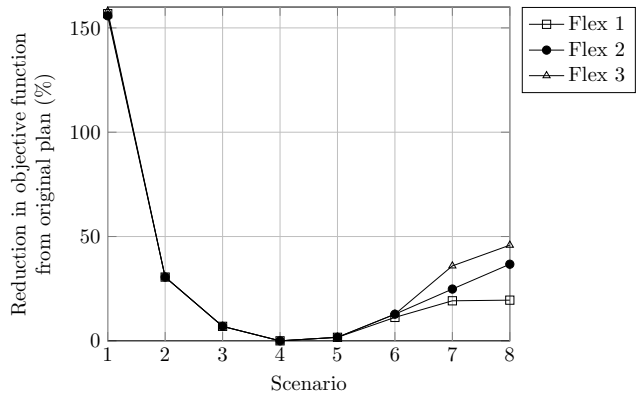
\begin{figure} [H]
\begin{center}
\begin{tikzpicture}[scale = 0.78]
\begin{axis}[
    ylabel style={align=center},
    ylabel =Reduction in objective function\\ from original plan (\%),
    ymin=0, ymax=160,
    xlabel={Scenario},
    xmin=1, xmax= 8,
    xtick=\empty,
    extra x ticks      ={1,2,3,4,5,6,7,8},
    extra x tick labels={1,2,3,4,5,6,7,8},
    legend pos = outer north east,
    ymajorgrids=true,
    grid, ]

    \addplot[mark=square] coordinates {(8,19.5)(7,19.19)(6,11.23)(5,1.69)(4,0)(3,6.95)(2,30.56)(1,156.8)};
            
    \addplot [mark=*] coordinates {(8,36.72)(7,24.82)(6,12.77)(5,1.69)(4,0)(3,6.94)(2,30.56)(1,155.8)};
        
    \addplot [mark=triangle] coordinates {(8,45.85)(7,35.92)(6,12.77)(5,1.69)(4,0)(3,6.94)(2,30.56)(1,158.02)};
    
    \legend{Flex 1, Flex 2,Flex 3},

\end{axis}
\end{tikzpicture}
\caption{Percentage reduction in objective function from original plan (No error in the beliefs) for different scenarios and flexibility levels.}
\label{reductRipening}
\end{center}
\end{figure}
%%%%%%%%%%%%%%%%%%%%%%%%%%%%%%%%%%%%%%%%%%%%%%%%%%%%%%%%%%%%
%\begin{tikzpicture}[scale = 0.78]
%\begin{axis}[
%    title={Objective function deficit},
%    ylabel={Objective function deficit (\$)},
%    ymin=0, ymax=55000000,
%    xlabel={1°stage decision error (\%)},
%    xmin=-85, xmax= 195,
%    xtick=\empty,
%    extra x ticks      ={-85,-55,-27,0,21,66,121,195},
%    extra x tick labels={-85,-55,-27,0,21,66,121,195},
%    legend pos = outer north east,
%    ymajorgrids=true,
%    grid, ]

%    \addplot[mark=square] coordinates 
%    {(195,21367500)(121,10378196)(66,4342044)(21,432253)
%    (0,0)(-27,1095991)(-55,4489330)(-85,12848444)};
            
%    \addplot [mark=*] coordinates 
%    {(195,50242750)(121,13422146)(66,4939447)(21,432253)
%    (0,0)(-27,1095698)(-55,4489330)(-85,13033944)};
        
%    \addplot [mark=triangle] coordinates 
%    {(195,50242750)(121,13422146)(66,4939534)(21,432253)
%    (0,0)(-27,1095698)(-55,4489330)(-85,13219444)};
    
%    \legend{Flex 1, Flex 2,Flex 3},

%\end{axis}
%\end{tikzpicture}

\begin{figure} [H]
\begin{center}
\begin{tikzpicture}[scale = 0.78]
\pgfplotsset{compat=1.9,compat/bar nodes=1.8,}
    \pgfplotstableread{
        Label series1 series2 
        1  44.18   55.82   
        2  31.15   68.85   
        3  25.33   74.67   
        4  00.00   00.00
        5  -22.24  -77.76   
        6  -20.64  -79.36   
        7 -18.15  -81.81    
        8 -07.29  -92.71    
    }\testdata
    \begin{axis}[
        xlabel={Scenario},
        ylabel style={align=center},
        ylabel=Percentage distribution of\\ increment/loss in income (\%),
        ybar stacked,
        ymin=-100,
        ymax=100,
        xtick=data,
        legend style={
            cells={anchor=west},
            legend pos= outer north east,
        },
        reverse legend=true,
        xticklabels from table={\testdata}{Label},
        xticklabel style={text width=2cm,align=center},
    ]
        \addplot [fill=black!80]
            table [y=series1, meta=Label, x expr=\coordindex]
                {\testdata};
                    \addlegendentry{Suboptimal harvest quality}
        \addplot [fill=gray!60]
            table [y=series2, meta=Label, x expr=\coordindex]
                {\testdata};
                    \addlegendentry{Unharvested grape}
    \end{axis}
\end{tikzpicture}
\caption{Distribution of income increment/loss between unharvested and sub-optimally harvested grapes for different scenarios.}
\label{yield_unreal2.2si}
\end{center}
\end{figure}
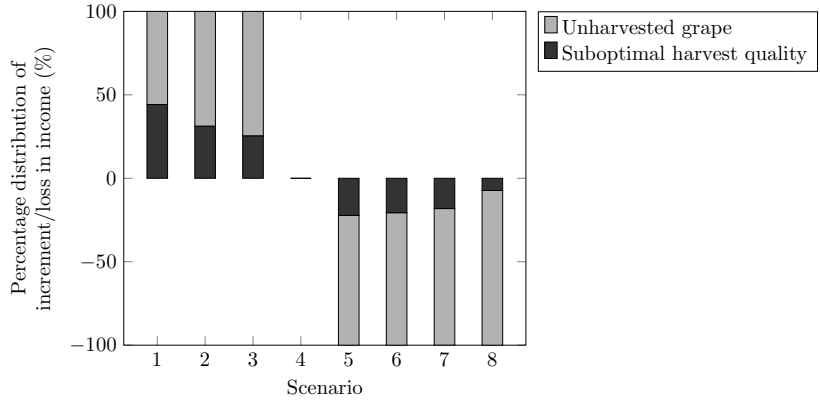

We can observe in figure \ref{yield_unreal2.2si} a similar behavior in comparison with the previous analysis. The main source of income losses when good scenarios probabilities are underestimated corresponds to the increment of unharvested grape. By the other hand, when good rate scenarios are overestimated, the compensating extra income is composed more equally between the increment of harvested grape and the better average quality of it. These results are almost the same for the other two levels of flexibility an so, we presented just the corresponding to flexibility level 2.
%%%%%%%%%%%%%%%%%%%%%%%%%%%%%%%%%%%%%%%%%%%%%%%%%%%%%%%%%%%%

These two last graphs show the results regarding to the grape performance separately by the ripening rates of grape, for flexibility two. High ripening rate correspond to grapes which its quality improves earlier, then medium ripening rate improves linearly and last, low ripening rate grape quality improves lately. Here, we can observe again that grapes which have a late improvement limit the decision-maker to adjust the harvest schedule according to its flexibility. We can notice this in low ripening rate blocks, where the unharvested amount and the loss income are higher than the others.
%%%%%%%%%%%%%%%%%%%%%%%%%%%%%%%%%%%%%%%%%%%%%%%%%%%%%%%%%%%%%%%%%%
\begin{figure} [H]
\begin{center}
\begin{tikzpicture}[scale = 0.78]
\begin{axis}[
    ylabel={Unharvested grape (\%)},
    ymin=-30, ymax=100,
    xlabel={Scenario},
    xmin=1, xmax= 8,
    xtick=\empty,
    extra x ticks      ={1,2,3,4,5,6,7,8},
    extra x tick labels={1,2,3,4,5,6,7,8},
    legend pos = outer north east,
    ymajorgrids=true,
    grid, ]

    \addplot[mark=square] coordinates {(8,50)(7,36.06)(6,21.73)(5,8.36)(4,0)(3,-11.03)(2,-28.63)(1,-19.76)};
            
    \addplot [mark=*] coordinates 
    {(8,16.07)(7,20.98)(6,13.43)(5,4.12)(4,0)(3,-7.65)(2,-21.16)(1,-10)};
    
    \addplot [mark=triangle] coordinates {(8,63.23)(7,40.4)(6,20.09)(5,6.4)(4,0)(3,-6.56)(2,-2.23)(1,-0.06)};
    
    \legend{High Ripening Rate, Medium Ripening Rate, Low Ripening Rate},

\end{axis}
\end{tikzpicture}
\caption{Percentage unharvested grape scenarios and ripening rates of grapes.}
\label{yield_unreal}
\end{center}
\end{figure}
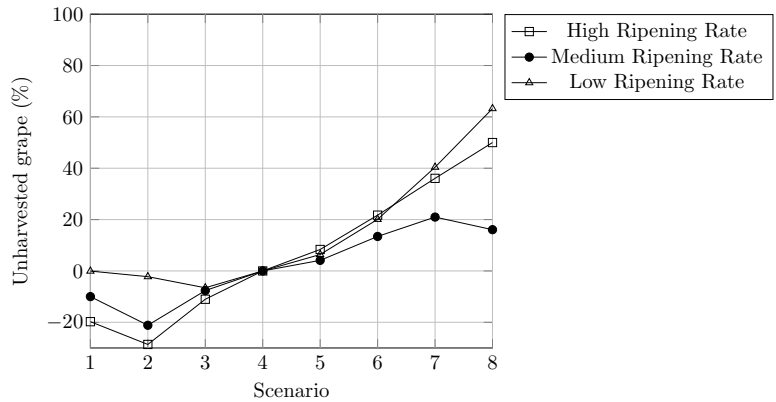

\begin{figure} [H]
\begin{center}
%%%%%%%%%%%%%%%%%%%%%%%%%%%%%%%%%%%%%%%%%%%%%%%%%%%%%%%%%%%%%%%%%%%%%%%%
\begin{tikzpicture}[scale = 0.78]
\begin{axis}[
    ylabel={Unrealized income (\$)},
    ymin=-1600000, ymax=35000000,
    xlabel={Scenario},
    xmin=1, xmax= 8,
    xtick=\empty,
    extra x ticks      ={1,2,3,4,5,6,7,8},
    extra x tick labels={1,2,3,4,5,6,7,8},
    legend pos = outer north east,
    ymajorgrids=true,
    grid, ]

    \addplot[mark=square] coordinates 
    {(8,12040000)(7,4934502)(6,2153388)(5,612888)(4,0)(3,-535716)(2,-1083566)(1,-399063)};
            
    \addplot [mark=*] coordinates 
    {(8,3870000)(7,4457135)(6,2243790)(5,503450)(4,0)(3,-674837)(2,-1586355)(1,-661042)};
    
    \addplot [mark=triangle] coordinates 
    {(8,31777000)(7,9224619)(6,4654959)(5,1056893)(4,0)(3,-660822)(2,-330308)(1,-10066)};
    
    \legend{High Ripening Rate, Medium Ripening Rate, Low Ripening Rate},

\end{axis}
\end{tikzpicture}
\caption{Unrealized income for different scenarios and ripening rates of grapes.}
\label{yield_unreal}
\end{center}
\end{figure}
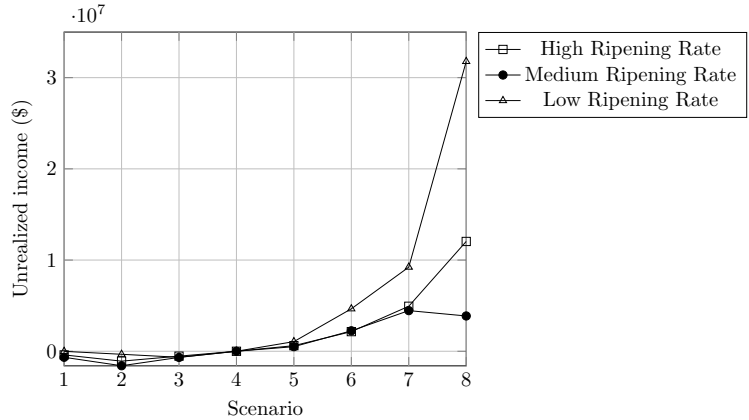

\section{Discussion and Conclusions}

In this research, we first develop a multistage optimization model, based on the one proposed by \cite{ferrer2008optimization} and differing from the one that \cite{avanzini2021comparing} present, since it considers the uncertainty in the yields and the transition probabilities. Using this model, we study the effect that errors in beliefs have on the value of using a Multistage Stochastic Optimization approach. Finally, we also analyze the effect that the quality and flexibility the resources have when errors in beliefs occur. From these findings, we can determine conditions under which acting to reduce the errors in future events will pay off, and conditions under which there is not much value to gain.

Results show that errors in grape-yield estimation have a significant impact on value and wine grape losses, evidencing that it is not symmetrical when yields are over- or under-estimated. Reductions in value can be up to 77.9\% and 232\% when the yields are over- or under-estimated by 100\%, respectively. They also show that an important part of these value losses comes from the reduction in net incomes, which can go up to 98\% of the potential amount. 

When yields are underestimated, flexibility or the ability to modify the decisions once the state of the nature reveals itself—for example, the process of adjusting the labor assignment—plays an important role. A high level of flexibility enables reducing the loss in value from the 77.9\% to just 36\%. 

This effect is mainly attributable to the reduction of net income losses, which can decrease from 98\% to 40\% as flexibility increases. This occurs because the loss in value is lightly compensated with an increment in the income, due to the available grapes. However, since there it is not enough labor to harvest the grapes, a significant amount is left unharvested and hence loss. If we look at how the model reacts in light of  quality, it generally favors the harvest of high-quality grapes over low quality, but the unrealized income is greater for the high-quality grapes, due to their greater value.

On the other hand, when yields are overestimated, results show that flexibility does not play a significant role in ameliorating the value loss. This is because the adjustment of labor cannot compensate for the reduction in the amount of available grapes to harvest. However, the harvest of grapes that would been unharvested under normal conditions, as well as the increment of high-quality grapes being harvested on their optimal dates due to labor availability, offset some of the value loss.

Looking at the different quality of grape in the contexts of under- and over-estimation of yields, the model tries to reduce comparatively the amount of unharvested high-quality grapes compared to those of low quality. Only when the yields decrease by 50\% is the percentage of high-quality grapes larger. If we compare value or unrealized gains from underestimated yields, most of the value loss comes from the high-quality grapes (due to their higher kg value) rather than those of low quality. However, in absolute terms, overestimated yields do not result in significant value loss.

For the same flexibility level, we analyzed the effect that the ripening rate has on the amount of unharvested grapes and on the unrealized gains when mistakes occur in the determination of the yields. When the ripening rate is low, grapes mature very close to their optimal date, and the yields are underestimated. The amount of unharvested grapes—hence, the unrealized gains—are largest, with insufficient labor to process them all. This amount decreases when the ripening rate is high with is a more ample window of time for harvesting the ripe grapes. This same effect occurs for the three levels of flexibility, and it increases as flexibility decreases.

We analyze the effect that errors in determining the transition probabilities have on the value of the plan. When the decision-maker is optimistic (the transition probabilities absorb the low yields), the value decreases significantly, with loss in value in the range of 50\% for the most optimistic scenario. When the decision-maker is pessimistic (the transition probabilities absorb the high yields), we see less reduction in value, compared to the optimistic scenario, which can go between 15\% and 30\% when we approach the extreme cases. 

If we look at the effect of flexibility on the value under different scenarios of transition-probability errors, results show that decision-maker optimism has no effect on reducing value loss. In pessimistic scenarios, flexibility ameliorates value loss by 15\%, due to workers’ ability to adjust productivity and reduce the amount of unharvested grapes.

For the case of the optimistic scenario, the reduction in value due to the change in the transition probabilities is attributable to high labor cost for the level of grape yields. Although reducing the unharvested grapes and harvesting them on dates closer to their optimal date produces some value, the reduction in yields and labor costs overcomes this value increment. In the pessimistic scenario, the value loss and unrealized gain results from the amount of unharvested grapes due to the lack of labor.

The maximum reachable quality of the grape also makes an important difference when errors in transition probabilities occur. When the decision-maker’s level of pessimism is high, the better average quality of the harvested grape lightly compensates for the value loss from the unharvested grape due to lack of labor. In these cases, the model prioritizes harvesting the high-quality grapes, leaving more low-quality grape unharvested. High-quality grapes still represent the greater part of harvest value.

If we analyze the effect that the ripening rates have on value when the scenario is optimistic, the excess available labor leaves less grape unharvested. When the decision-maker is too pessimistic, grapes with a low ripening rate are more difficult to harvest, with the smaller window available for harvesting on optimal dates preventing their harvest starting earlier.

\subsection{Managerial Insights}

According to our results, a manager facing the possibility that his beliefs can be incorrect or inaccurate would likely hire labor in excess of what planning recommended. This would allow him to face a potential underestimation of yields without incurring high-cost hiring at the last minute,  forcing the leaving of unharvested grape due to lack of labor. This same behavior appears in other industries (e.g., airlines, where adding robustness to the flight crew pairs solutions by also adding ground time between the incoming and outgoing flights) \citep{ehrgott2002constructing, yen2006stochastic, weide2010iterative, dunbar2012robust}.

Errors regarding yield estimation seem to have a larger effect on value than making the same mistakes in transition probabilities. Thus, facing the decision of where to invest in producing more reliable information, the decision should aim toward accuracy in determining grape yields.

Flexibility adds significant value to the planning process, as \cite{avanzini2021comparing} indicate. So, if possible, the manager should try to induce flexibility into the plan, making possible the adjusting of the plan as the states of nature appear. This could occur by directly relating the payments of the workers to productivity and reducing hiring and lay-off costs. Finally, achieving a higher level of flexibility can occur if the decision-maker has a buffer for harvesting capacity by hiring extra workers. 

Regarding the maximum achievable quality of the grapes, the harvesting plan should prioritize higher-quality grapes since these determine most of the income and net utilities of production. Also, a situation with a lack of available labor calls for adapting the harvest plan to leave as little high-quality grape unharvested as possible. Small increments of unharvested premium grape can cause higher income losses than low-quality grape produces. 

In the case of ripening rates, the preferred possibility is to plant varieties that have high ripening rates and wider optimal harvest windows. This increases the flexibility of the decision-maker to adjust the harvest plan by starting the harvesting earlier than expected when reduced labor is available. Likewise, avoiding varieties with low ripening rates maintains the decision-maker’s flexibility to modify the harvest plan.

Finally, the results this work obtains show that the effects of errors in beliefs on the planning value behave asymmetrically. This indicates that using regular optimization techniques to solve the problem at hand is not the best method. A better mechanism to handle errors in beliefs could be the use of a robust optimization approach \citep{ben2009robust} or a conditional value at risk \citep{rockafellar2002conditional}, which minimizes losses under worst-case scenarios or average expected costs considering them all, a more accurate approach to results more robust against uncertainty.

\subsection{Future Research}

This research focuses on analyzing the errors regarding yields and transition probability matrix beliefs. The agricultural realm includes many other errors in beliefs, such as: labor productivity levels, quality level of grapes, probability of adverse climatic events (e.g., rain). Further research that studies the effects of these errors would be a great contribution to the winery industry and agriculture in general.

Transition probabilities and scenarios were considered as discrete events. This work could go even further if continuous distribution functions were associated with the yields. Also, the same could apply to the change of quality through periods, or ripening rates behavior. This would enable implementing a distributionally robust optimization model approach, and a closer representation to the real distribution could offer more useful managerial insights. 

Our results indicate that the use of a labor buffer to act as a reserve, in cases of underestimating the harvested stock, is a desirable feature. The use of a robust optimization approach or a conditional value at risk helps by defining optimal buffer levels.

Finally, studying the effects of using different types of machine-learning schemes in improving the precision and accuracy of updating would reduce the gap between the believed values and the real ones. Since the quality and amount of information directly relates to the quality of the forecasts, studying the use and value of sensor applications on agricultural industries to collect better and larger amounts of information to support decision-making would be desirable.

%\section{Acknowledgments}
%Alejandro F. Mac Cawley acknowledges financial support from FONDECYT Iniciacion Project \# 11180502. \\

\newpage

\bibliography{references}

\end{document}